\documentclass{siamart}

\usepackage{lipsum}
\usepackage{amsfonts}
\usepackage{graphicx}
\usepackage{amsopn}
\usepackage{currfile}
\usepackage{enumerate}

\usepackage{color}
\usepackage{colortbl} % color for tables
\usepackage{multirow}

\usepackage{hyperref}

\usepackage{changes}
\definecolor{mygray}{RGB}{30,10,110}
\newcommand{\tc}{\textcolor{black}}

\usepackage{epstopdf}
\usepackage{algorithmic}

\newsiamremark{remark}{Remark}

\newcommand{\N}{\mathbb N}

\newcommand{\R}{\mathbb R}

\newcommand{\mcM}{\mathcal M}

\newcommand{\rank}{\operatorname{rank}}
\newcommand{\tr}{\operatorname{tr}}
\newcommand{\sym}{\operatorname{sym}}
\newcommand{\Skew}{\operatorname{skew}} %\skew exists already in LaTeX
\newcommand{\dist}{\operatorname{dist}}
\newcommand{\Exp}{\operatorname{Exp}}
\newcommand{\Log}{\operatorname{Log}}
\newcommand{\Cay}{\operatorname{Cay}}

\makeatletter
\renewcommand*\env@matrix[1][*\c@MaxMatrixCols c]{%
  \hskip -\arraycolsep
  \let\@ifnextchar\new@ifnextchar
  \array{#1}}
\makeatother

\newcount\Comments  % 0 suppresses notes to selves in text
\newcommand{\kibitz}[2]{\ifnum\Comments=1\textcolor{#1}{#2}\fi}

\newcommand{\e}{\epsilon}

\ifpdf
  \DeclareGraphicsExtensions{.eps,.pdf,.png,.jpg}
\else
  \DeclareGraphicsExtensions{.eps}
\fi

\newcommand{\TheTitle}{Computing the Riemannian logarithm on the Stiefel manifold: metrics, methods and performance} 
\newcommand{\TheAuthors}{Ralf Zimmermann and Knut H{\"u}per}

\headers{Computing the Riemannian Stiefel log}{\TheAuthors}

\title{{\TheTitle}} %\thanks{This work was funded by the XXX under contract no.~XX.}}

\author{
  Ralf Zimmermann\thanks{Department of Mathematics and Computer Science, University of Southern Denmark (SDU), Odense, Denmark,
    (\email{zimmermann@imada.sdu.dk})}
 \and 
 Knut H{\"u}per\thanks{Institute of Mathematics,
 Julius-Maximilians-Universit{\"a}t, W{\"u}rzburg, Germany,
 (\email{hueper@mathematik.uni-wuerzburg.de})}
}

\ifpdf
\hypersetup{
  pdftitle={\TheTitle},
  pdfauthor={\TheAuthors}
}
\fi

\begin{document}

\maketitle

% REQUIRED
\begin{abstract}
We address the problem of computing Riemannian normal coordinates on the real, compact Stiefel manifold of orthonormal frames. The Riemannian normal coordinates are based on the so-called Riemannian exponential and the associated Riemannian logarithm map and enable to transfer almost any computational procedure to the realm of the Stiefel manifold.
To compute the Riemannian logarithm is to solve the (local) geodesic endpoint problem.
Instead of restricting the consideration to geodesics with respect to a single selected metric, we consider a family of Riemannian metrics introduced by H\"uper, Markina and Silva-Leite that includes the Euclidean and the canonical metric as prominent examples.

As main contributions, we provide (1) a unified, structured, reduced formula for the Stiefel geodesics.
\tc{The formula is unified in the sense that it works for the full family of metrics under consideration. It is structured in the sense that it relies on matrix exponentials of skew-symmetric matrices exclusively. It is reduced in relation to the dimension of the matrices of which matrix exponentials have to be calculated.
We provide} (2) a unified method to tackle the geodesic endpoint problem numerically and
(3) we improve the existing Riemannian log \tc{algorithm} under the canonical metric \tc{in terms of the computational efficiency}.
The findings are illustrated by means of numerical examples, where the novel algorithms prove to be the most efficient methods known to this date.
\end{abstract}

% REQUIRED

% REQUIRED

\section{Introduction}
\label{sec:intro}
\tc{Riemannian computing methods have established themselves as important tools in a large variety of applications, including  computer vision, machine learning, and optimization, see
\cite{AbsilMahonySepulchre2004,AbsilMahonySepulchre2008,BegelforWerman2006,EdelmanAriasSmith1999,Gallivan_etal2003,  Lui2012, Rahman_etal2005, Rentmeesters2013} and the anthologies \cite{Minh:2016:AAR:3029338, RiemannInComputerVision}.
They also gain increasing attention in statistics and data science
\cite{patrangenaru2015nonparametric} %, trendafilov2021multivariate}
and in numerical methods for differential equations \cite{BennerGugercinWillcox2015,Celledoni_2020,hairer06gni, iserles_munthe-kaas_norsett_zanna_2000,Zimmermann_MORHB2021}.\\
One way to enable the practical execution of data processing methods on a curved manifold $\mcM$ is via
working in local coordinates.}
This holds among others for basic tasks like averaging, clustering, interpolation and optimization.
Of special importance are the {\em Riemannian normal coordinates}, as they are {\em radially isometric} \cite[\S 5]{Lee2018riemannian}.
The Riemannian normal coordinates rely on the {\em Riemannian exponential} and the {\em Riemannian logarithm}, which are local diffeomorphisms:
The exponential at a manifold location $p\in\mcM$ sends a tangent vector $v$ (i.e., the velocity vector of a manifold curve) to the endpoint $q=c(1)$ of a geodesic curve $c$ that starts from $p=c(0)$ with velocity $v=\dot c(0)$.
The Riemannian logarithm at $p$ maps a manifold location $q\in\mcM$ to the starting velocity vector $v$ of a geodesic $c$ that connects $p=c(0)$ and $q=c(1)$.
\Cref{fig:ExpLogStiefel} illustrates this process.
Hence, the Riemannian logarithm is associated with the {\em geodesic endpoint problem}:
\begin{quote}
\tc{``Given $p,q\in\mcM$, 
find a geodesic arc that connects $p$ and $q$.''}
\end{quote}
In this work, we tackle the local geodesic endpoint problem on the Stiefel manifold of orthonormal frames. Geodesics depend on the way \tc{the length of} velocity vectors of curves are measured and thus on the Riemannian metric. Popular choices for metrics on the Stiefel manifold are the {\em Euclidean metric} and the {\em canonical metric}. These will be detailed in \Cref{sec:Stiefel}.
Rather than restricting the considerations to either of these two,
we work with the one-parameter family of metrics that is introduced
in \cite{HueperMarkinaLeite2020}. This family contains the Euclidean and the canonical Stiefel metric as special cases.
\paragraph{Original contributions}
\begin{itemize}
\item We start from the results of \cite{EdelmanAriasSmith1999} and             
      \cite{HueperUllrich2018} and derive a \tc{unified formula} for the Stiefel geodesics and thus for the Stiefel exponential. 
      \tc{Here, unified is to be understood in the sense that the formula} works for all metrics in the one-parameter family under consideration.
      Moreover it features the same \tc{skew-symmetric} structure as exhibited by the canonical geodesics and comes at roughly the same computational costs.
\item We provide new theoretical insights on the structure of the matrices that come 
      into consideration as candidate solutions for the Stiefel geodesic endpoint problem. For rectangular matrices of dimensions $(n\times p)$, $n\gg p$, this reduces the endpoint problem to finding suitable $(p\times p)$-orthogonal matrices.
\item We provide efficient algorithms for computing the Riemannian logarithm in a 
      unified way for the full one-parameter family of metrics of \cite{HueperMarkinaLeite2020}.  
\item For the special case of the canonical metric, the endpoint problem features a simplified structure that was exploited in \cite{StiefelLog_Zimmermann2017}. We refine  this approach and thus accelerate the existing method.
\item We juxtapose the various methods by means of numerical experiments, \tc{where the new methods prove to outperform both their predecessors from \cite{Bryner2017,StiefelLog_Zimmermann2017} as well as Newton-based approaches}.
\end{itemize}
%accelerated computational schemes when compared to both 
%\cite{StiefelLog_Zimmermann2017}, \cite{Bryner2017}.

\paragraph{Related work}
% Literature on computing the Riemannian logarithm on the Stiefel manifold is sparse.
The reference \cite[section 5]{Rentmeesters2013} tackles the local
geodesic endpoint problem for the canonical metric via a Riemannian optimization approach;
\cite{StiefelLog_Zimmermann2017} also works in the setting of the canonical metric
and provides a matrix-algebraic algorithm based on the Baker-Campbell-Hausdorff formula with guaranteed local linear convergence.
An algorithm for computing the Stiefel logarithm for the Euclidean metric
is considered in \cite{Bryner2017} and is based on the general ``shooting method'', see \cite[section 6.5]{SrivastavaKlassen2016}.
\tc{The thesis \cite{Sutti_PhD_2020} considers single-shooting and multiple-shooting methods based on Newton's method to solve the geodesic endpoint problem under the canonical metric. It also features an algorithm for computing global Stiefel geodesics that is
based on the ``leapfrog method'' of \cite{noakes_1998} for general manifolds, see also the associated preprint \cite{Sutti_V:2020b}.} This approach requires methods \tc{that tackle the geodesic endpoint problem locally, i.e., for input points that are close enough to each other, as building blocks. The notion of being ``close enough'' depends on geometric quantities (like the injectivity radius) but also on the numerical algorithm that is applied to the problem.}
\paragraph{Notational specifics}
\label{sec:Notation}
For $p\in \N$, the $(p\times p)$-identity matrix is denoted by $I_p\in\R^{p\times p}$, or simply $I$, if the dimension is clear.
The $(p\times p)$-orthogonal group, i.e., the set of all
square orthogonal matrices is denoted by
\[
  O(p) = \{\phi \in \R^{p\times p}| \phi^T\phi = \phi\phi^T = I_p\}.
\]
The standard matrix exponential and matrix logarithm are denoted by
\[
 \exp_m(X):=\sum_{j=0}^\infty{\frac{X^j}{j!}}, \quad \log_m(I+X):=\sum_{j=1}^\infty{(-1)^{j+1}\frac{X^j}{j}}.
\]
% We use the symbols
% $\Exp, \Log$ for the Riemannian counterparts on the Stiefel manifold.
The sets of symmetric and skew-symmetric $(p\times p)$-matrices are $\sym(p) = \{A\in\R^{p\times p}|A^T=A\}$ and
$\Skew(p) = \{A\in\R^{p\times p}|A^T=-A\}$, respectively.
Overloading this notation, $\sym(A) = \frac12(A+\tc{A^T}), \Skew(A) = \frac12(A-A^T)$
denote the symmetric and skew-symmetric parts of a matrix $A$.

Unless stated otherwise, when we employ the QR-decomposition of a rectangular matrix $A\in\R^{n\times p}$, 
we implicitly assume that $n\geq p$ and work with
the \tc{`compact'} QR-decomposition $A=QR$, with $Q\in \R^{n\times p}$, $R\in \R^{p\times p}$.
%
%
%
%
%%
%%%
%%%%
%%%%%
%%%%%%
%%%%%%%
\section{The Stiefel manifold}
\label{sec:Stiefel}
This section reviews the essential aspects of Stiefel manifolds
in regards of numerical, matrix-algorithmic applications.
For additional background, see \cite{AbsilMahonySepulchre2008,EdelmanAriasSmith1999,Zimmermann_MORHB2021}.
For featured applications, see, e.g., \cite{bernstein2012tangent, Chakraborty2017, HueperRollingStiefel2008, Turaga_2008}.

The {\em Stiefel manifold} $St(n,p)$ 
is the set of rectangular, column-orthonormal $n$-by-$p$ matrices,
\[
  St(n,p):= \{U \in \R^{n\times p}| \hspace{0.1cm}  U^TU = I_p\}, \quad p\leq n.
\]
Observe that this matrix set is the pre-image $St(n,p)=F^{-1}(0)$ of the function
$F:\R^{n\times p}\rightarrow \sym(p), Y\mapsto Y^TY-I_p$.
By the regular value theorem, it is a differentiable manifold of dimension $np-\frac12 p(p+1)$.

% An alternative approach is via Example \ref{ex:QuotBasic}, where $St(n,p)$ is identified with the quotient space
% $St(n,p)\cong O(n)/(I_p\times O(n-p))$
% under actions of the closed subgroup $I_p\times O(n-p):= \left\{\begin{pmatrix}
%                             I_p &0\\
%                              0    & R
%                            \end{pmatrix}|\hspace{0.1cm} R\in O(n-p)\right\}\leq O(n)$.
% Two square orthogonal matrices in $O(n)$ are identified 
% as the same point on $St(n,p)$, if their first $p$ columns coincide, see \cite[\S 2.4]{EdelmanAriasSmith1999}.
%

The {\em tangent space} $T_USt(n,p)$ at $U \in St(n,p)$ is represented as
\[
  T_USt(n,p) = \left\{\Delta \in \R^{n\times p} |\quad U^T\Delta \in\Skew(p)\right\}.
\]
For brevity, we will often write $T_U$ instead of $T_USt(n,p)$.
Every tangent vector $\Delta \in T_U$ may be written as
\begin{align}
\label{eq:tang_St1}
   \Delta&= UA + (I-UU^T)T,& A \in  \Skew(p), & \quad T\in\R^{n\times p} \mbox{ arbitrary,}\\
\label{eq:tang_St2}  
   \Delta&= UA + U^{\bot}H,& A \in \Skew(p), & \quad H\in\R^{(n-p)\times p} \mbox{ arbitrary.}
\end{align}
In the latter case, $U^{\bot}\in St(n,n-p)$ is an orthonormal completion such that the matrix $\begin{pmatrix}[c|c]U & U^{\bot}\end{pmatrix}$
is orthogonal.
% Since different $T\neq\widetilde{T}$ may produce the same projection $(I-UU^T)T=(I-UU^T)\widetilde{T}$,
% only \eqref{eq:tang_St2} reflects the actual degrees of freedom (DoFs) in $T_U$,
% namely $\frac{1}{2}p(p-1)$ (DoFs in $A\in\Skew(p)$) plus $(n-p)p$ (DoFs in $H\in \R^{(n-p)\times p}$).
%
Any matrix $W\in \R^{n\times p}$ can be projected onto $T_U$ by
\begin{equation}
 \label{eq:tangproj}
 \Pi_U(W) = W - U\sym(U^TW),
\end{equation}
see \cite[eqs. (2.3), (2.4)]{EdelmanAriasSmith1999}.

In order to turn $St(n,p)$ into a {\em Riemannian manifold} \cite{Lee2018riemannian},
a metric, i.e., an inner product $\langle\cdot,\cdot \rangle_U$ 
with associated norm $\|\cdot \|_U = \sqrt{\langle\cdot,\cdot \rangle_U}$ on the tangent spaces $T_U$ must be defined for all $U\in St(n,p)$.
Given a metric, the length of a curve $C:[a,b]\to St(n,p)$ is
$ L(C):= \int_a^b \|\dot C(t) \|_{C(t)}dt.$
Candidates for length-minimizing curves are called {\em geodesics}
and are locally
% stationary points of the length functional
% associated with a regular variation of a family of curves with fixed endpoints. Geodesics $t\mapsto C(t)$ are 
uniquely determined by an
ordinary initial value problem when specifying a starting point $C(0)$ and a starting velocity $\dot C(t)$, \cite[\S 6]{Lee2018riemannian}.
It is obvious that geodesics depend on the underlying Riemannian metric.
Geodesics give rise to the Riemannian exponential map.
On the Stiefel manifold, the Riemannian exponential at a base point $U\in St(n,p)$ sends a Stiefel tangent vector $\Delta$
to the endpoint $C(1) = \widetilde{U}\in St(n,p)$ of a geodesic $t\mapsto C_{U,\Delta}(t)$ that starts from $C(0) = U$ with velocity vector $\dot C(0) = \Delta$,
\[
 \Exp_U(\Delta) := C_{U,\Delta}(1). %, \quad C_{U,\Delta}(t) = \Exp_U(t\Delta).
\]
\tc{As a consequence, $\Exp_U(t\Delta) = C_{U,\Delta}(t)$ for $t\in[0,1]$}. Knowing the Riemannian exponential is knowing the geodesics and vice versa.

The Riemannian exponential is locally invertible.
The inverse is called the {\em Riemannian logarithm} and is denoted
\begin{equation}
 \label{eq:RiemannLog}
 \Log_{U}: St(n,p)\ni\widetilde{U}\mapsto \Log_{U}(\widetilde U):= \Exp_U^{-1}(\widetilde U) \in T_USt(n,p).
\end{equation}
More precisely, $\Log_U(\widetilde U)$ is well defined for all 
$\widetilde U$ within the {\em injectivity radius} $r_{St}(U)$ of $St(n,p)$ at $U$.
The injectivity radius at $U$ is the Riemannian distance of $U$ to its cut locus $C_U$. The cut locus, in turn, is the set of points beyond which the geodesics starting from $U$ cease to be length-minimizing
% . 
% The global injectivity radius is $r_{St} := \inf \{r_{St}(U)|\hspace{0.1cm} U\in St(n,p)\}=\inf_{U\in St(n,p)} \{\dist(U, C_U)\}$,
 \cite[p. 271]{DoCarmo2013riemannian}.
%
%---------------------------------------------------------------------------
\begin{figure}[ht]
\centering
\includegraphics[width=0.9\textwidth]{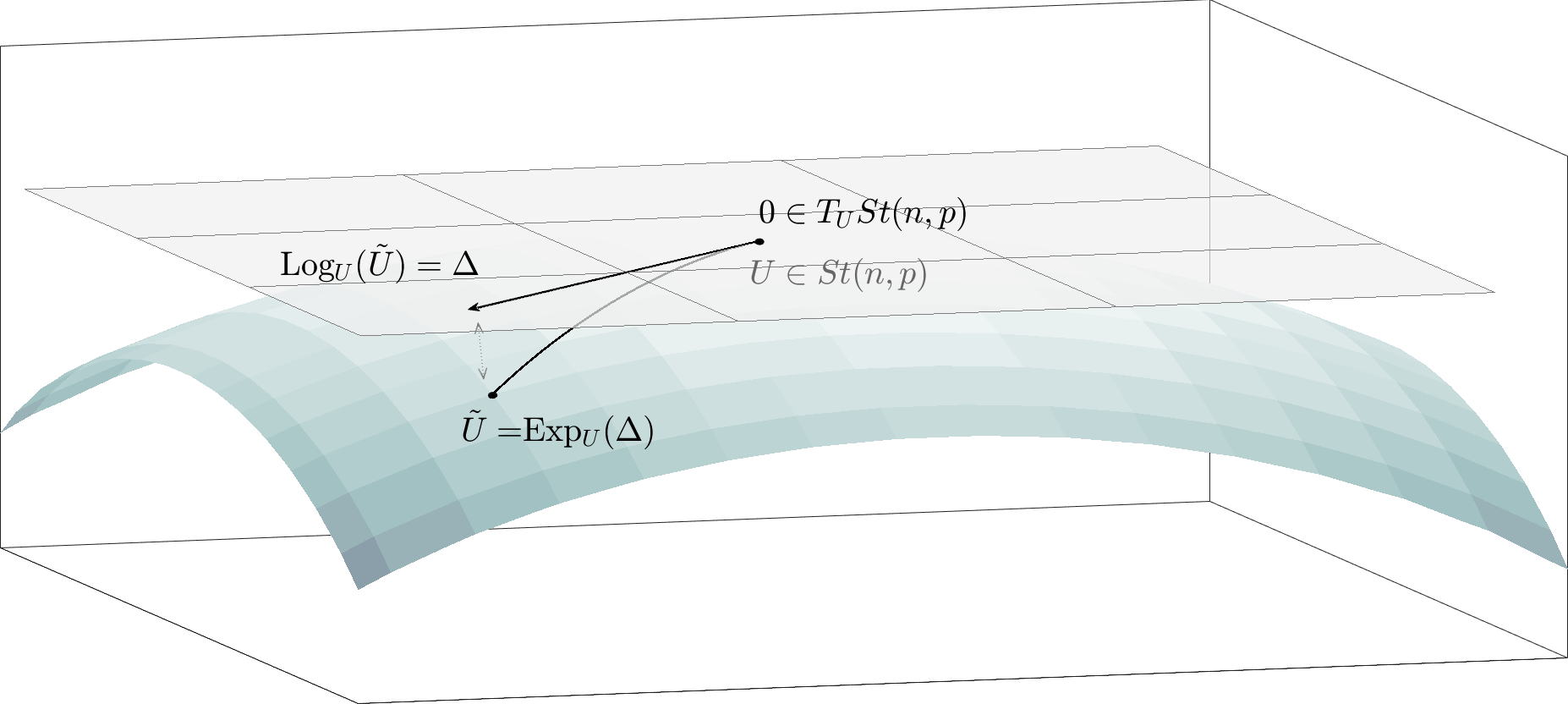}
\caption{(cf. \Cref{sec:Stiefel}) Conceptual visualization of the Stiefel manifold $St(n,p)$ (curved surface) and associated tangent space $T_USt(n,p)$ (shaded plane).
The Riemannian exponential sends a tangent vector $\Delta$ from $T_USt(n,p)$ to the endpoint $\widetilde U$ of the geodesic $C_{U,\Delta}:[0,1]\to  St(n,p)$
with initial values $C_{U,\Delta}(0) = U, \dot C_{U,\Delta}(0) = \Delta$. 
The Riemannian logarithm inverts this process.
}
\label{fig:ExpLogStiefel}
\end{figure}
%------------------------------------------------------
Combined, the Riemannian logarithm and exponential provide the {\em Riemannian normal coordinates}, which allow to map data points back and forth between the curved 
manifold and the flat tangent space, see \Cref{fig:ExpLogStiefel}. This is crucial for all data processing 
operations on manifolds (optimization, interpolation, averaging, clustering,...).
The Riemannian normal coordinates are special in that they are length-preserving along geodesic rays; one speaks of radial isometries.
\subsection*{The Euclidean and the canonical metric and generalizations}
\label{sec:Stiefel_metrics}
Let $U\in St(n,p)$ and let $\Delta= UA + U^{\bot}H$, $\widetilde{\Delta}= U\widetilde{A} + U^{\bot}\widetilde{H} \in T_U$. Here and in the following, $A = U^T\Delta, \widetilde A = U^T\widetilde\Delta \in \Skew(p)$.
There are two standard metrics on the Stiefel manifold.\\
The {\em Euclidean metric} on $T_U$ is the one inherited from the ambient $\mathbb{R}^{n\times p}$:
 \[
   \langle\Delta, \widetilde{\Delta}\rangle_U^e = \tr( \Delta^T \widetilde{\Delta}) = \tr A^T\widetilde{A} + \tr H^T\widetilde{H}.
 \]
The {\em canonical metric} on $T_U$  is derived from the
 quotient representation $St(n,p) = O(n)/(O(n-p))$ of the Stiefel manifold, see \cite{EdelmanAriasSmith1999}, and reads
  \[
  \langle \Delta, \widetilde{\Delta}\rangle_U^c = \tr\left(\Delta^T(I-\frac{1}{2}UU^T)\widetilde{\Delta}\right) 
  = \frac12 \tr A^T\widetilde{A} +  \tr H^T\widetilde{H}.
 \]
\tc{Let $A = (a_{ij})_{i,j\leq p}$ and $H = (h_{ij})_{i\leq n; j\leq p}$.}
The Euclidean metric corresponds to measuring
tangent vectors $\Delta=UA+U^\bot H$ in the Frobenius matrix norm 
\[
\sqrt{\langle \Delta,\Delta \rangle_U^e}= \|\Delta\|_F = \sqrt{\|A\|_F^2 + \|H\|_F^2}
= \sqrt{2\sum_{i<j} a_{ij}^2 + \sum_{i,j} h_{ij}^2},
\]
while the canonical metric yields
\[
\sqrt{\langle \Delta,\Delta \rangle_U^c}= %\sqrt{\|\Delta\|_F^2 - \frac12\|A\|_F^2}=
\sqrt{\frac12\|A\|_F^2 + \|H\|_F^2}
 = \sqrt{\sum_{i<j} a_{ij}^2 + \sum_{i,j} h_{ij}^2}.
\]
In this sense, the Euclidean metric disregards the skew-symmetry of $A$ 
and \tc{the independent entries $a_{ij}, i <j$ are counted} twice, as was observed in \cite[\S 2.4]{EdelmanAriasSmith1999}.

%
%%
%%%
%%%% alpha metrics
The work \cite{HueperMarkinaLeite2020} recognizes the Euclidean and the canonical metric
as special cases of a one-parameter family of inner products
\begin{equation}
  \label{eq:alphaMetrics}
  \langle \Delta, \widetilde{\Delta}\rangle_U^\alpha = \tr\left(\Delta^T(I-\frac{2\alpha+1}{2(\alpha+1)}UU^T)\widetilde{\Delta}\right) 
  = \frac{1}{2(\alpha+1)} \tr A^T\widetilde{A} +  \tr H^T\widetilde{H},
   \end{equation}
for $\alpha\in \R\setminus\{-1\}$.\footnote{For a deeper reason, why $\alpha=-1$ has to be excluded, see \cite{HueperMarkinaLeite2020}.}
For $\alpha = -\frac12$ and $\alpha = 0$,
the Euclidean and the canonical metric are recovered, respectively.
\tc{As can be seen from \eqref{eq:alphaMetrics}, the metric parameter $\alpha$ controls how much weight is put on the $\frac12(p-1)p$ degrees of freedom in the matrix $A\in \Skew(p)$ relative to the $(n-p)p$ degrees of freedom in $H\in \R^{(n-p)\times p}$.
Perfect balance is at $\alpha = 0$, which corresponds to the canonical metric.}
%
%%
%%%
\subsection*{Calculating the Riemannian Stiefel exponential: The state of the art}
\label{sec:Stiefel_exp}
A closed-form expression for the Stiefel exponential w.r.t. the Euclidean metric is derived in 
\cite[\S 2.2.2]{EdelmanAriasSmith1999},
  \begin{equation}
  \label{eq:StExpEAS}
   \widetilde{U} = \Exp_U^e(\Delta) = 
   \begin{pmatrix}U & \Delta \end{pmatrix}
   \exp_m
     \begin{pmatrix}
   A & -\Delta^T\Delta\\
   I_p & A
  \end{pmatrix}
     \begin{pmatrix}
   I_p\\
   0
  \end{pmatrix}
  \exp_m(-A).
  \end{equation}
In \cite{HueperUllrich2018}, an alternative formula is obtained:
  \begin{equation}
  \label{eq:StExpHU}
   \widetilde{U} = \Exp_U^e(\Delta) = 
   \exp_m\left(
   \Delta U^T - U\Delta^T
   \right)
   U
  \exp_m(-A).
  \end{equation}
The advantage is that in this form, the Stiefel exponential features
only matrix exponentials of skew-symmetric matrices.
The downside is that $\Delta U^T - U\Delta^T\in \Skew(n)$
and working with $(n\times n)$-matrices might be prohibitively expensive
in large-scale applications, where $n\gg p$.

In \cite{HueperMarkinaLeite2020}, the formula \eqref{eq:StExpHU} is generalized to the 
exponential for all $\alpha$-metrics of \eqref{eq:alphaMetrics},
  \begin{equation}
  \label{eq:StExpAlpha}
   \widetilde{U} = \Exp_U^\alpha(\Delta) = 
   \exp_m\left(-\frac{2\alpha +1}{\alpha +1}UA U^T + \Delta U^T - U\Delta^T \right)U
  \exp_m\left(\frac{\alpha}{\alpha +1}A\right).
  \end{equation}

An algorithm for computing the Stiefel exponential w.r.t. the canonical metric
was derived in \cite[\S 2.4.2]{EdelmanAriasSmith1999}:
Given $U, \Delta$, first compute a compact QR-decomposition
$QR = (I-UU^T)\Delta$ with $Q\in St(n,p), R\in \R^{p\times p}$. Then form
  \begin{equation}
  \label{eq:StExpCanEAS}
   \widetilde{U} = \Exp_U^c(\Delta) = 
   \begin{pmatrix}U & Q \end{pmatrix}
   \exp_m
     \begin{pmatrix}
   A & -R^T\\
   R & 0
  \end{pmatrix}
     \begin{pmatrix}
   I_p\\
   0
  \end{pmatrix},
  \end{equation}
where again $A=U^T\Delta \in \Skew(p)$.
Using this form becomes efficient, if $p<\frac{n}{2}$.
%
%
%
%%
%%%
\subsection*{Calculating the Riemannian Stiefel logarithm: State of the art}
\label{sec:Stiefel_log}
\begin{algorithm}[t]
\caption{Shooting method, adapted from \cite[Alg. 1]{Bryner2017}}
\label{alg:SimpleShooting}
\begin{algorithmic}[1]
  \REQUIRE{Stiefel matrices $U, \widetilde{U}\in St(n,p)$, convergence threshold $\e>0$,
  time step array $T=\{t_0, t_1,\ldots, t_m\}$, with $t_0=0, t_m = 1.0$.}% (discretized unit interval).}
  \STATE{ $\gamma \leftarrow  \|\widetilde U -U\|$} \hfill \COMMENT{\tc{compute gap between base point and target}}
  \STATE{ $\Delta \leftarrow \gamma \frac{\Pi_U(\widetilde U)}{\|\Pi_U(\widetilde U)\|}$}
    \hfill \COMMENT{project gap vector $\widetilde U -U$ onto $T_U$, preserve length}
  \WHILE{ $\gamma> \e$}
    \FOR{$j=1,\ldots,m$}
    \STATE{ $\widetilde U^s(j) \leftarrow \Exp_U(t_j\Delta)$}\hfill \COMMENT{discrete representation of geodesic}
    \ENDFOR
    %\STATE{$\widetilde U^s \leftarrow\Exp_U(\Delta)$}  \hfill \COMMENT{hit of current shot}
    \STATE{ $\Delta^s \leftarrow \widetilde U^s(m) - \widetilde U $}
            \hfill \COMMENT{current gap vector: compare endpoint of geodesic to $\widetilde U$}
    \STATE{$\gamma \leftarrow\|\Delta^s\|$}
    \FOR{$j=m,\ldots,0$}
    \STATE{ $\Delta^s\leftarrow \gamma \frac{\Pi_{\widetilde U^s(j)}(\Delta^s)}{\|\Pi_{\widetilde U^s(j)}(\Delta^s)\|}$}
    \hfill \COMMENT{initial projection plus approximate parallel transport}
    \ENDFOR
    \STATE{ update $\Delta\leftarrow\Delta - \Delta^s$}
  \ENDWHILE
  \ENSURE{$\Delta$}

 \textbf{Note:} The first inner iteration in the loop in steps 9-11 projects the gap vector of step 7 onto $T_{\widetilde U^s(m)}$.
 The remaining inner iterations in this loop realize an approximation of the parallel transport along the discretized geodesic $t\mapsto \Exp_U(t\Delta)$.
\end{algorithmic}
\end{algorithm}
Computing the Riemannian logarithm corresponds to solving the {\em geodesic endpoint problem} locally:
\begin{quote}
 Given $U, \widetilde U\in St(n,p)$ find a starting velocity $\Delta\in T_U$ such that
 \[
  \Exp_U(\Delta) = \widetilde U \quad \left(\Leftrightarrow  \Delta = \Log_{U}(\widetilde U)\right).
 \]
\end{quote}
A generic method for solving the geodesic endpoint problem 
on any manifold is the shooting method, see \cite[section 6.5]{SrivastavaKlassen2016}.
This method is considered in \cite{Bryner2017} to compute the Stiefel logarithm. For the specific case of $St(4,2)$ it features in \cite{Sundaramoorthi2011}.
The generic principle of the shooting method is as follows:
\begin{itemize}
 \item find an initial guess $\Delta_0\in T_U$.
 \item ``shoot'' a geodesic in the direction of $\Delta_0$, i.e., compute $\widetilde{U}_0 := \Exp_U(\Delta_0)$.
 \item measure the gap  $f_g$ between $\widetilde{U}_0$ and the actual target $\widetilde U$ as a function of $\Delta$.
 \item use information on the gap to update $\Delta_1 \leftarrow \Delta_0$ and repeat.
\end{itemize}
There are many options to implement these steps.
In \cite{Bryner2017}, the initial guess is chosen as $\Delta_0 = \Pi_U(U-\widetilde U)\in T_U$. The gap is measured in the ambient space $\R^{n\times p}$ as
$f_g(\Delta)  = \|\widetilde{U}_0-\widetilde U\|_F^2 = \|\Exp_U(\Delta)-\widetilde U\|_F^2$.
A natural choice for updating the shooting direction $\Delta$ is a gradient descent based on this gap function.
\tc{An alternative is to tackle the matrix root finding problem $\Exp_U(\Delta)-\widetilde U = 0$ with the Newton method.
This approach is pursued in \cite[Section 2.3]{Sutti_PhD_2020}.
However, both of these approaches require} in particular the derivative of the matrix exponential,
%$D(\exp_m)_X(Z)$ 
which is expensive to obtain, see \cite[Section 10.6]{Higham:2008:FM} 
for the general formulas and \cite[Prop.12]{Sundaramoorthi2011}, \cite[Lem. 5]{ZimmermannHermite_2020} for precise applications
to the Stiefel exponential.

The reference \cite{Bryner2017} proposes a method that avoids calculating derivatives. The idea is to project the gap vector 
$\Exp_U(\Delta)-\widetilde U$ from $\R^{n\times p}$ onto $T_{\widetilde U}St(n,p)$ and to parallel-translate the result along the geodesic $t\to\Exp_U(t\Delta)$ back to $T_{U}St(n,p)$, \tc{where the update of the shooting direction $\Delta$ is then performed.
This process is detailed in \cref{alg:SimpleShooting}.
For a higher computational efficiency, only an approximation of the parallel transport is realized
on equidistant time steps $0=t_0, t_1,\ldots, t_m=1$ in the unit interval $[0,1]$.
The more time steps there are, the more accurate will be the result of the parallel transport.
However, this does not help, if the tangent information that is to be transported is of poor quality in the first place.\\
The ``multiple shooting'' as featured in \cite[Section 2.4]{Sutti_PhD_2020}
shares the idea of dissecting the geodesic lines under consideration into $m$ segments.
Yet in this method, one actually solves subproblems on the segments with Newton's method.}
%
%
% 
% If one would `shoot' from a point $v$ at a point $\widetilde v$ in a Euclidean space, in some direction $d$, then the first shot
% would lead to $\widetilde v^s = v +d$ with a gap of $d^s = \widetilde v^s - \widetilde v$. Then, the update would be corrected to
% $d \leftarrow d - d^s$ which directly gives the correct shooting direction. \cref{alg:SimpleShooting} mocks this principle and iterates.
%
%

A Stiefel log algorithm that is tailored for the canonical metric 
and features guaranteed local linear  convergence is given in \cite{StiefelLog_Zimmermann2017}.
This approach will be explained and enhanced in \Cref{sec:AlgImpro}.
\section{Fast computational schemes for solving the geodesic endpoint problem}
\label{sec:FastComp}
In this section, we first prepare the grounds for a unified framework to tackle the local geodesic endpoint problem \tc{on the Stiefel manifold} under the one-parameter family of metrics. Then, we introduce practical numerical algorithms for 
computing the \tc{associated Stiefel logarithm}, which improve on the state of the art in terms of the computational efficiency.
\subsection{A reduced formula for the Stiefel exponential}
\label{sec:reduceExpEucl}
As a starting point for efficient computational schemes, we derive an alternative
expression for the  $\alpha$-metric Stiefel exponential that combines the advantages of \cref{eq:StExpEAS} and \cref{eq:StExpAlpha} and represents a considerable
computational reduction if $n\gg p$.
To this end, assume that $p\leq \frac{n}2$. Let $U\in St(n,p)$ be given and take any suitable column-orthonormal extension $U^\bot\in St(n-p,p)$, i.e., any $U^\bot$ such that
$\Phi = \begin{pmatrix}[c|c]
  U & U^\bot
 \end{pmatrix}\in O(n)$
is square and orthogonal.
Observe that $\Phi^TU = 
\begin{pmatrix}
 I_p\\ 0
\end{pmatrix}
$.
With $\nu = \frac{2\alpha +1}{\alpha +1}$, $\mu = \frac{\alpha}{\alpha +1}= \nu -1$, we rewrite \cref{eq:StExpAlpha}:
\begin{eqnarray}
\nonumber
 \Exp_U^\alpha(\Delta) 
 &=& \Phi\Phi^T
 \exp_m\left(-\nu UA U^T + \Delta U^T - U\Delta^T \right)
 \Phi\Phi^TU
  \exp_m(\mu A)
\\
\nonumber
  &=& \Phi
   \exp_m\left(
   -\nu \Phi^TUA U^T \Phi+
   \Phi^T(\Delta U^T - U\Delta^T)\Phi
   \right)\Phi^T
   U
  \exp_m(\mu A)\\
\nonumber
  &=&  \Phi
   \exp_m
   \begin{pmatrix}
    -\nu A + U^T\Delta - \Delta^T U & - \Delta^T U^\bot\\
    (U^\bot)^T\Delta & 0
   \end{pmatrix}
  \begin{pmatrix}
    I_p\\
    0
  \end{pmatrix}
  \exp_m(\mu A)\\
  &=&
  \begin{pmatrix}[c|c]
  U & U^\bot
 \end{pmatrix}
 \exp_m
   \begin{pmatrix}
    (2-\nu)A & -H^T\\
    H & 0
   \end{pmatrix}
  \begin{pmatrix}
   I\\
    0
  \end{pmatrix} \exp_m(\mu A)
  , \quad H = (U^\bot)^T\Delta.\label{eq:curious1}
\end{eqnarray}
Let 
$H = \widetilde Q \begin{pmatrix}B\\0 \end{pmatrix}$
with $\widetilde Q= \begin{pmatrix}[c|c]\widetilde Q_p & \widetilde Q_{n-2p}\end{pmatrix} \in O(n-p)$ be
a (full) QR-decomposition \tc{of $H$}.
(Only the orthonormality of $\widetilde Q$ matters, $B$ need not be triangular for the following considerations.)
We can factor
\[
    \begin{pmatrix}
    (2-\nu)A & -H^T\\
    H & 0
   \end{pmatrix}
   =
      \begin{pmatrix}
    I_p & 0\\
    0 & \widetilde Q
   \end{pmatrix}
   \exp_m
      \begin{pmatrix}
    (2-\nu)A & -B^T & 0\\
    B  & 0    & 0\\
    0  & 0    & 0
   \end{pmatrix}
      \begin{pmatrix}
    I_p & 0\\
    0 & \widetilde Q^T
   \end{pmatrix},
\]
which yields 
$
 \Exp_U^\alpha(\Delta) = 
 \begin{pmatrix}[c|c]
  U & U^\bot \widetilde Q_p
 \end{pmatrix}
  \exp_m
   \begin{pmatrix}
    (2-\nu)A & -B^T\\
    B & 0
   \end{pmatrix}
  \begin{pmatrix}
   \exp_m(\mu A)\\
    0
  \end{pmatrix} 
$.
Note that $U^\bot H = U^\bot(U^\bot)^T\Delta = (I-UU^T)\Delta$.
Hence, instead of computing a QR-de\-com\-po\-si\-tion of $U^\bot H$,
we can directly compute a compact QR-decomposition  $(I-UU^T)\Delta = QB$ with $Q\in St(n,p)$
(Again, no special structure is required for $B\in \R^{p\times p}$.)
This leads to the following proposition.
\begin{proposition}
\label{prop:StExpAlpha}
Let \tc{$\alpha \neq -1$}. 
 For $U\in St(n,p)$, $\Delta\in T_U St(n,p)$ the Stiefel exponential \tc{under the $\alpha$-metric of \eqref{eq:alphaMetrics} }reads
 \begin{equation}
  \label{eq:StAlphaRed}
  \Exp_U^\alpha(\Delta) = 
 \begin{pmatrix}[c|c]
  U& Q
 \end{pmatrix}
  \exp_m
   \begin{pmatrix}
    \frac{1}{\alpha+1}A & -B^T\\
    B & 0
   \end{pmatrix}
  \begin{pmatrix}
   I_p\\
    0
  \end{pmatrix}\exp_m\left(\frac{\alpha}{\alpha+1} A\right),
 \end{equation}
 where $A = U^T\Delta\in \Skew(p)$ and $QB= (I-UU^T)\Delta\in \R^{n\times p}$
 is any matrix decomposition with 
 with $Q\in St(n,p)$ and $B\in \R^{p\times p}$.
\end{proposition}
As with \eqref{eq:StExpCanEAS}, this form features only standard matrix exponentials of skew-symmetric 
matrices of size $p$ rather than $n$.
With \Cref{prop:StExpAlpha}, the Stiefel exponential is computable for
all metrics in the $\alpha$-family \cref{eq:alphaMetrics} in $\mathcal{O}(np^2)$ flops. 
A more detailed look on the calculations reveals that the formula \cref{eq:StAlphaRed}
remains valid also if $p>\frac{n}2$. However, in this case, it represents in fact an increase in dimension for the main matrix exponential rather than a reduction when compared 
to \eqref{eq:StExpAlpha}.
The formula \cref{eq:StAlphaRed} continues to hold if the orthogonal component $(I-UU^T)\Delta$ of the tangent vector is rank-deficient.
In this case, it is understood that the QR-decomposition
$QB= (I-UU^T)\Delta$ be arranged such that 
$QB = \begin{pmatrix}[c|c]Q_r & Q_{p-r}\end{pmatrix}
\begin{pmatrix}B_r \\ 0\end{pmatrix}$, where
$B_r\in \R^{r\times p}$ and $\rank(B_r) = \rank((I-UU^T)\Delta)=r$ and
$B_r = 0$ in the extreme case of $0 = (I-UU^T)\Delta$.

Before we continue with the main body of this work,
we note an interesting aside.
\begin{lemma}
 \label{lem:curious}
 For $A\in \Skew(p)$, $H\in \R^{(n-p)\times p}$, it holds
 \[
  \begin{pmatrix}
   I_p & A\\
   0 & H
  \end{pmatrix}
  \exp_m\begin{pmatrix}
         A & A^2 - H^TH\\
         I_p & A
        \end{pmatrix}
  \begin{pmatrix}
   I_p\\
   0
  \end{pmatrix}
  = \exp_m
  \begin{pmatrix}
   2A & -H^T\\
   H & 0
  \end{pmatrix}
  \begin{pmatrix}
   I_p\\
   0
  \end{pmatrix}.
 \]
\end{lemma}
\begin{proof}
 Fix an orthogonal matrix $\begin{pmatrix}
                            U & U^\bot
                           \end{pmatrix}\in O(n)$ and construct $\Delta = UA + U^\bot H$.
 Then, the lemma is a consequence of the fact that both \eqref{eq:StExpEAS} and \eqref{eq:curious1} with $\alpha = -\frac12$ are valid expressions of 
 the (unique) Stiefel matrix exponential $\Exp^e_U(\Delta)$.
 \tc{The underlying geodesics $t\mapsto \Exp^e_U(t\Delta)$
 satisfy the same initial value problem and thus coincide, see \cite[Theorem 4.27, p. 103]{Lee2018riemannian}.}
\end{proof}
% To try and prove \Cref{lem:curious} directly seems to be difficult.
%
%%
%%%
%%%%
%%%%%
%%%%%%
%%%%%%%
%%%%%%%%
\subsection{The Stiefel log matrix equations}
\label{sec:matrixEq}
We continue to work under the assumption that $p\leq \frac{n}2$ with the 
setting of $n\gg p$ in mind for practical big-data applications.
In this section, we show that, in essence, computing the Stiefel logarithm for all $\alpha$-metrics corresponds to solving a nonlinear matrix equation,
where the unknown is a $2p\times 2p$ skew-symmetric matrix.

Let $U\in St(n,p)$ and let $\widetilde U \in St(n,p)$ be within the injectivity radius 
$r_{St}(U)$ of $St(n,p)$ at $U$.
Let $\Delta \in T_U$ be such that $\Exp_U^\alpha(\Delta) = \widetilde U$.
Then, $\widetilde U$ has a representation
\[ \widetilde U = UM + QN,\quad M,N\in \R^{n\times p}, \quad M^TM + N^TN = I_p,
\]
where $Q\in St(n,p)$ is the `$Q$-factor' of a compact QR-decomposition $QB=(I-UU^T)\Delta$.
For the canonical metric ($\alpha=0$), this is immediately clear from \eqref{eq:StExpCanEAS}.
For the full family of $\alpha$-metrics, this follows from \eqref{eq:StAlphaRed}.
More precisely, $M,N$ form the first block column of an orthogonal matrix
\begin{equation}
\label{eq:OrthMatpExp}
\begin{pmatrix}
  M & X\\
  N & Y
 \end{pmatrix}
    = 
\exp_m\begin{pmatrix}
         \frac{1}{\alpha +1} A &  -B^T\\
         B &  0
 \end{pmatrix}
 \begin{pmatrix}
         \exp_m(  \frac{\alpha}{\alpha +1} A) &  0\\
         0 &  I_p
 \end{pmatrix} \in SO(2p).
\end{equation}
Recall $A = U^T\Delta\in \Skew(p)$ and that 
$\alpha=0$ and $\alpha = -\frac{1}{2}$ reproduce the canonical and the Euclidean metric, 
respectively. For brevity, write again
\[\mu = \frac{\alpha}{\alpha +1}, \quad \alpha \neq -1.\]

Now, we aim at \tc{finding the inverse  $(\Exp_U^\alpha)^{-1}(\widetilde U)=\Delta$} with $\Delta$ as unknown.
We can still obtain $M$ and $N$ from the given data points $U,\widetilde U$ via
\[
 M = U^T\widetilde U, \quad QN = (I-UU^T)\widetilde U.
\]
For the latter equation, any matrix decomposition that 
represents $(I-UU^T)\widetilde U$ via a `subspace factor' $Q$ with orthonormal columns and a corresponding `coordinates factor' $N$
is suitable in order to obey the constraint that $M$ and $N$ be blocks of an orthogonal matrix, i.e., $I = M^TM + N^TN$.
This incorporates some ambiguity.
If $(I-UU^T)\widetilde U$ has full rank $p$, then $Q$ and $N$ are unique up to 
a rotation/reflection $\Phi\in O(p)$. More precisely, for $\widetilde Q \widetilde N =(I-UU^T)\widetilde U = QN$,
it holds $Q = \widetilde Q\Phi, N = \Phi^T \widetilde N$ with $\Phi = \widetilde Q^TQ$.
If $\rank (I-UU^T)\widetilde U= r <p$, then we may assume that 
\begin{equation}
 (I-UU^T)\widetilde U = QN =
 \begin{pmatrix}[c|c]
  Q_r &Q_{p-r}
 \end{pmatrix}
 \begin{pmatrix}
  N_r\\
  0
 \end{pmatrix},
\end{equation}
where $Q_r \in St(n,r), Q_{p-r}\in St(n,p-r), N_r \in \R^{p\times r}$.
This determines $Q_r, N_r$ uniquely up to 
a rotation/reflection $\Phi_r\in O(r)$. Yet, the column block $Q_{p-r}$ may be an arbitrary orthonormal extension.

Once $Q, M,N$ are chosen and fixed, we can compute matrix blocks $X_0$, $Y_0$
such that they form an orthogonal completion 
$V=\begin{pmatrix}
  M & X_0\\
  N & Y_0
 \end{pmatrix}\in O(2p)$.
The restriction of the exponential map to the skew-symmetric matrices
\[
 \exp_m|_{\Skew(p)}: \Skew(p)\rightarrow  SO(p)
\]
is surjective \cite[\S. 3.11, Thm. 9]{godement2017introduction}.
Hence, \eqref{eq:OrthMatpExp} requires to select $X_0, Y_0$ such that $\det(V) = +1$.\footnote{As an aside, this is why the example S1
from the supplements of \cite{StiefelLog_Zimmermann2017} fails.}
Again, such a completion is not unique and it cannot be expected that the such found $X_0, Y_0$ align with the structure of \eqref{eq:OrthMatpExp}.
Yet, there is an orientation preserving orthogonal matrix $\Phi = \exp_m(C), C\in\Skew(p)$ such that
$\begin{pmatrix}
  X_0\\
  Y_0
 \end{pmatrix}\Phi
 =
 \begin{pmatrix}
  X\\
  Y
 \end{pmatrix}
 $
is the sought-after completion.
In summary, computing the Stiefel logarithm boils down to solving the following nonlinear matrix equation:

\begin{eqnarray}
\label{eq:OrthMatpLog1}
\text{solve} &&
F\begin{pmatrix}
  A & -B^T\\
  B & C
 \end{pmatrix}
 = 
\begin{pmatrix}
  M & X_0\\
  N & Y_0
 \end{pmatrix},
 \quad \text{ for }
 F: \Skew(2p) \to SO(2p),\\
 \nonumber %\label{eq:OrthMatpLog2}
  F\begin{pmatrix}
  A & -B^T\\
  B & C
 \end{pmatrix}
&  = &
\exp_m\begin{pmatrix}
         (1-\mu) A &  -B^T\\
         B &  0
 \end{pmatrix}
 \begin{pmatrix}
         \exp_m(\mu A) &  0\\
         0 &  \exp_m( -C)
 \end{pmatrix}, \mu = \frac{\alpha}{\alpha +1}.
\end{eqnarray}
The sought-after tangent vector is then
$\Delta = UA + QB\in T_U$.
It is worth mentioning that any ambiguity in $Q$ has no impact on the final $\Delta$, not even in the case where the component $(I-UU^T)\widetilde U$ of $\widetilde U$ 
does not feature full rank.
This is confirmed in the next theorem, which can be considered as a generalization of 
\cite[Thm. 3.1]{StiefelLog_Zimmermann2017}.
\begin{theorem}
\label{thm:LogSearchSpace}
Let $U,\widetilde U\in St(n,p)$ with $\dist(U,\widetilde U) < r_{St}(U)$ and let 
\[
U^T\widetilde U = M ,
 \quad  (I-UU^T)\widetilde U =QN
 =
 \begin{pmatrix}
  Q_r & Q_{p-r}
 \end{pmatrix}
 \begin{pmatrix}
  N_r\\
  0
 \end{pmatrix}, \quad r\leq p,
\]
with $Q\in St(n,p)$ and $r=\rank(N_r) = \rank(N)$.
Then $\Delta = \Log_U^\alpha(\widetilde U)$ features a representation
\[
 \Delta =
 UA + 
 Q
 B
 = UA + \begin{pmatrix}
  Q_r& Q_{p-r}
 \end{pmatrix}
  \begin{pmatrix}
  B_r\\
  0
 \end{pmatrix}, \quad A\in \Skew(p),  \quad B_r \in \R^{r\times p}.
\]
Here, $A,B$ are components of a skew-symmetric block matrix that solves \eqref{eq:OrthMatpLog1}.
\end{theorem}
\begin{proof}
 Let \tc{$M, N \in \R^{p\times p}$, with $N = \begin{pmatrix}
  N_r\\
  0
 \end{pmatrix}, N_r\in \R^{r\times p}$} of full rank $r$ constructed as stated above.
 We can restrict the considerations to block extensions $X_0,Y_0$ to an orthogonal matrix of the form
 \[
V = \begin{pmatrix}
  M   & X_0\\
  N   & Y_0
 \end{pmatrix}
 =\begin{pmatrix}[c|c|c]
  M   & X_r^0 & 0\\ \hline
  N_r & Y_r^0 & 0\\ \hline
  0   & 0   & I_{p-r}
 \end{pmatrix}
 \in SO(2p),%\quad N_r\in \R^{r\times p}, X_r^0 \in \R^{p\times r}, Y_r^0 \in \R^{r\times r}.
 \]
  where  $X_r^0 \in \R^{p\times r}, Y_r^0 \in \R^{r\times r}$ are obtained, say, via the Gram-Schmidt method.\\
 Let 
 $\begin{pmatrix}
  A & -B^T\\
  B & C
 \end{pmatrix}$
 be a solution to \eqref{eq:OrthMatpLog1} that preserves the above structure of $V$
 so that in particular 
 $
 \begin{pmatrix}[c|c]
  X_r^0 & 0\\ \hline
  Y_r^0 & 0\\ \hline
  0   & I_{p-r}
 \end{pmatrix}
 \exp_m(C)
 = 
  \begin{pmatrix}[c|c]
  X_r & 0\\ \hline
  Y_r & 0\\ \hline
  0   & I_{p-r}
 \end{pmatrix}
 $. This entails 
 $
 \exp_m(C) = \Phi = \begin{pmatrix}
  \Phi_r & 0\\
  0 & I_{p-r}
 \end{pmatrix}\in SO(p)
 $.
 Then, for all $\alpha$-metrics,
 it holds with $\mu =\frac{\alpha}{\alpha+1}$ that
 \begin{equation}
 \label{eq:OrthMatpLog2}
 %  \log_m\left(
%   V
%   \begin{pmatrix}
%   \Psi & 0\\
%   0    & \Phi
%  \end{pmatrix}
%  \right)
%  %
%  = 
 \log_m\left(
 \begin{pmatrix}[c|c|c]
  M   & X_r^0 & 0\\ \hline
  N_r & Y_r^0 & 0\\ \hline
  0   & 0   & I_{p-r}
 \end{pmatrix}
 \begin{pmatrix}[c|c|c]
   \exp_m(-\mu A)& 0      & 0\\ \hline
  0   & \Phi_r & 0\\ \hline
  0   &    0   & I_{p-r}
 \end{pmatrix} \right)
 =
 \begin{pmatrix}
   (1-\mu)A & -B^T\\
  B & 0
 \end{pmatrix}.
 \end{equation}
Because the matrix logarithm acts block-wise on block-diagonal matrices, the block structure of the matrix on the left hand side entails a corresponding
block structure for the right-hand side of \eqref{eq:OrthMatpLog2}.
In particular, necessarily 
$B =   \begin{pmatrix}
  B_r\\
  0
 \end{pmatrix}$ with $\rank(B_r) = r$.
 
 Now, define $\Delta = UA + Q_r B_r$ and \tc{evaluate} the Stiefel exponential to check if
 $\Exp_U^\alpha(\Delta) = \widetilde U$.
 According to \eqref{eq:StAlphaRed}, the procedure 
 requires to compute $U^T\Delta = A$, as well as an orthogonal decomposition
 \[
 (I-UU^T)\Delta  
 =
 \begin{pmatrix}[c|c]
   \widehat Q_r & \widehat Q_{p-r} 
  \end{pmatrix}
  \begin{pmatrix}
  \widehat B_r\\
  0
 \end{pmatrix}.
 \]
 The matrix factors $\widehat Q_r$ and $\widehat B_r$ are not necessarily equal to $Q_r$ and $B_r$, respectively. Yet, by construction, it holds $\widehat Q_r \widehat B_r=(I-UU^T)\Delta  =(I-UU^T) Q_r B_r = Q_r B_r$.
 Since $Q_r, \widehat Q_r$ span the same column-space,
 $S:= Q_r^T\widehat Q_r$ is orthogonal and $\widehat Q_r = Q_r S, \widehat B_r = S^T B_r$.
 The Stiefel exponential produces
 \begin{eqnarray*}
 \Exp_U^\alpha(\Delta) &=&
   \begin{pmatrix}[c|c|c]
   U & \widehat Q_r & \widehat Q_{p-r} 
  \end{pmatrix}
  \exp_m
   \begin{pmatrix}
  (1-\mu) A & -\widehat B^T_ r & 0\\
  \widehat B_r & 0 & 0\\
  0&0&0
 \end{pmatrix}
  \begin{pmatrix}
         \exp_m(\mu A)\\
         0\\
         0
 \end{pmatrix}\\
 &=& 
    \begin{pmatrix}[c|c]
   U & \widehat Q_r
  \end{pmatrix}
  \exp_m
   \begin{pmatrix}
   (1-\mu) A& -\widehat B^T_ r\\
  \widehat B_r & 0
 \end{pmatrix}
  \begin{pmatrix}
         \exp_m(\mu A)\\
         0
 \end{pmatrix}\\
 &=&
   \begin{pmatrix}[c|c]
   U & \widehat Q_rS
  \end{pmatrix}
  \exp_m
   \begin{pmatrix}
  (1-\mu) A & -\widehat B^T_ r S\\
  S^T \widehat B_r & 0
 \end{pmatrix}
  \begin{pmatrix}
         \exp_m(\mu A)\\
         0
 \end{pmatrix}\\
  &=&
   \begin{pmatrix}[c|c]
   U &  Q_r
  \end{pmatrix}
  \exp_m
   \begin{pmatrix}
  (1-\mu)A & - B^T_ r\\
   B_r & 0
 \end{pmatrix}
  \begin{pmatrix}
         \exp_m(\mu A)\\
         0
 \end{pmatrix}
 = \begin{pmatrix}[c|c]
   U &  Q_r
  \end{pmatrix}
    \begin{pmatrix}
         M\\
         N_r
 \end{pmatrix} =  \widetilde U
 \end{eqnarray*}
according to \eqref{eq:OrthMatpLog2}.
\end{proof}
\begin{remark}[The search space for the Stiefel log algorithm]
\label{rem:searchSpaceStLog}
For Stiefel points $U$ and $\widetilde{U} = \Exp_U^\alpha(\Delta)$ within the injectivity radius at $U$, \Cref{thm:LogSearchSpace} shows that both the location $\widetilde U$ and the tangent vector $\Delta$
are in the same matrix space 
\begin{equation}
\label{eq:LogSearchSpace}
 \mathcal{S} = \{U V| V\in \R^{p\times p}\}\oplus \{Q W| W\in \R^{p\times p}\}\subset \R^{n\times p}.
\end{equation}
This fact can be exploited in numerical schemes for computing the Stiefel logarithm.
A direct consequence is that any numerical algorithm can focus on finding the missing factors $A,B\in \R^{p\times p}$.
\end{remark}
%
%%
%%%
%%%%
%%%%%
%%%%%%
%%%%%%%
%%%%%%%%
\subsection{A {\em p}-shooting method tailored to the Stiefel logarithm}
\label{sec:accShootStLog}
In this section, we customize the generic iterative shooting method of \Cref{alg:SimpleShooting} for solving the geodesic endpoint problem.
With large ``tall-and-skinny'' matrices and big data applications in mind, our contribution is to modify the required calculations such that the loop iterations
are tailored for the use of \eqref{eq:StAlphaRed} and \eqref{eq:LogSearchSpace} and work exclusively  with $2p\times 2p$ matrices rather than with $n\times p$ matrices.
This leads to considerable savings, if $n\gg p$.

Given $U,\widetilde U\in St(n,p)$, the starting point is the representation 
 $\widetilde U = UU^T\widetilde U + (I-UU^T)\widetilde U = U\widehat M + Q\widehat N$, where
 $\widehat M =U^T\widetilde U\in \R^{p\times p}$ and $Q\widehat N= (I-UU^T)\widetilde U \in \R^{n\times p}$ is a compact QR-decomposition.
The essential observation is that all iterates $\Delta, \Delta^s$ of \cref{alg:SimpleShooting} actually remain in the matrix space $\mathcal{S}$ of \eqref{eq:LogSearchSpace} that is 
spanned by the fixed $U$ and $Q$.
According to \Cref{thm:LogSearchSpace} and \Cref{rem:searchSpaceStLog}, $\mathcal{S}$ also contains the sought-after solution.

\begin{proposition}
 \label{prop:stay_p}
 Let $U,\widetilde U\in St(n,p)$ and let and $Q\widehat N= (I-UU^T)\widetilde U \in \R^{n\times p}$ be a compact QR-decomposition.
 All tangent matrices $\Delta$ produced iteratively by \cref{alg:SimpleShooting} and all the update 
 matrices $\Delta^s$ that are the final outcome of one pass through the while-loop in \cref{alg:SimpleShooting}  allow for a representation of the form
 \[
  \Delta = UA + QR, \quad \Delta^s = UA^s + QR^s, \quad  A,A^s\in \Skew(p), R, R^s\in \R^{p\times p}.
 \]
\end{proposition}
{\bf Remark:} It is important to emphasize that it is {\em the same} Q-factor that works for all the iterates $\Delta, \Delta^s$. This enables to restrict the update to the $A$- and $R$-factors.
The intermediate matrices $\Delta^s$ under the while-loop are also contained in the matrix space $\mathcal{S}$ of \Cref{rem:searchSpaceStLog}, but they are of the form $\Delta^s = UX^s + QR^s$ with $X^s \not\in \Skew(p)$ in general.

\begin{proof}
 The initial $\Delta_0$ is obtained from the projection
 \begin{eqnarray*}
  \Delta_0 &=& \Pi_U(\widetilde U - U) = \Pi_U(\widetilde U) 
         = U\widehat M + Q\widehat N - U\sym(\widehat M)\\
         &=& U\Skew(\widehat M) + Q\widehat N =: UA_0 + Q R_0. 
 \end{eqnarray*}
 The first update $\Delta^s$ is obtained by projecting the gap vector $\widetilde U^s_0 - U$ onto $T_{\widetilde U}$ and then further onto $T_U$ by following the geodesic `backwards in time', see steps 9--11 in \cref{alg:SimpleShooting}.
 At every time instant $t$, it holds
 \[
  \widetilde U^s(t) = \Exp_U^\alpha(t\Delta_0) = \begin{pmatrix}[c|c]
                                     U & Q
                                    \end{pmatrix}
                                     \begin{pmatrix}
                                     M_t &  X_t\\
                                     N_t &  Y_t
                                    \end{pmatrix}
                                     \begin{pmatrix}
                                     I_p \\
                                     0
                                    \end{pmatrix}
                                    = UM_t + QN_t,
 \]
 where $\begin{pmatrix}
         M_t &  X_t\\
         N_t &  Y_t
        \end{pmatrix}$
is evaluated according to \eqref{eq:OrthMatpExp}
with inputs $A = tA_0$, $B = tR_0$.

If any matrix of the form $W =UX+QY$ is projected onto the tangent space at any $UM_t + QN_t$,
the result is
\begin{eqnarray*}
 \Pi_{UM_t + QN_t}(W) 
 &=& U(X - M_t\sym(M_t^T X + N_t^TY)) + Q(Y-N_t\sym(M_t^T X + N_t^TY))\\
 &=& UX^\Pi +  QY^\Pi.
\end{eqnarray*}
Hence, all operations performed in \cref{alg:SimpleShooting} take place in the matrix space
$\mathcal{S}  = \{UX + QY| X,Y\in \R^{p\times p}\}$ of \eqref{eq:LogSearchSpace}.
%
%
% \begin{algorithm}
% \caption{Exp4Geo}
% \label{alg:Exp4Geo}
% \begin{algorithmic}[1]
%   \REQUIRE{$A\in \Skew(p)$, $B\in \R^{p\times p}$}
%   %
%   \IF{metric == 'canonical'}
%     \STATE{T}
%   \ELSIF{metric == 'Euclidean'}
%     \STATE{T}
%   \ENDIF
%   %
%   \STATE{}
%   \ENSURE{}
% %
% \end{algorithmic}
% \end{algorithm}
%
%
%
% A calculation shows
%  \begin{eqnarray*}
%   \Delta^s_0 &=& \Pi_U(\Pi_{\widetilde U}(\widetilde U^s - \widetilde U)) = \Pi_U(\Pi_{\widetilde U}(\widetilde U^s))\\
%   &=& \Pi_U\left((U(M_0-\widehat M\sym(\widehat M^TM_0+\widehat N^TN_0)) + Q(N_0-\widehat N\sym(\widehat M^TM_0+\widehat N^TN_0))\right)\\
%   &=& U\Skew(M_0-\widehat M\sym(\widehat M^TM_0+\widehat N^TN_0)) + Q(N_0-\widehat N\sym(\widehat M^TM_0+\widehat N^TN_0))\\
%   &=:& UA^s_0 + QR^s_0. 
%  \end{eqnarray*}
% This gives the next iterate
% \[
%  \Delta_1 = \Delta_0 - \Delta^s_0 = U(A_0-A^s_0) + Q(R_0 - R^s_0)=: UA_1 + QR_1.
% \]
% The proof is completed by a straightforward induction.
\end{proof}

With \cref{prop:stay_p}, the computational costs associated with the shooting method can be reduced considerably.
This gives rise to \cref{alg:p-Shooting}.
\begin{algorithm}
\caption{p-Shooting method}
\label{alg:p-Shooting}
\begin{algorithmic}[1]
  \REQUIRE{Stiefel matrices $U, \widetilde{U}\in St(n,p)$, convergence threshold $\e>0$,
  array $T=[t_0, t_1,\ldots, t_m], t_0=0, t_m = 1.0$ (discretized unit interval), metric parameter $\alpha$.}
  \STATE{$\widehat M = U^T\widetilde U$}
  \STATE{$Q\widehat N = \widetilde U - U\widehat M$}\hfill \COMMENT{compact QR-decomposition}
  \STATE{ $\gamma \leftarrow  \sqrt{\|\widehat M-I_p\|^2 + \|\widehat N\|^2}$}\hfill \COMMENT{note $\widetilde U -U = U(\widehat M -I_p)+ Q\widehat N$}
  \STATE{ $A \leftarrow \frac{\gamma\Skew(\widehat M)}{\sqrt{\|\Skew(\widehat M)\|^2 + \|\widehat{N}\|^2}}$,
  $\quad   R \leftarrow\frac{ \gamma\widehat N}{\sqrt{\|\Skew(\widehat M)\|^2 + \|\widehat{N}\|^2}}$,}
  \STATE{ }
  \WHILE{ $\gamma> \e$}
    \FOR{$j=1,\ldots,m$}
    \STATE{ \tc{ $\begin{pmatrix}
              M(t_j)\\ N(t_j) 
             \end{pmatrix} \leftarrow
             \exp_m\left(t_j\begin{pmatrix}
                \frac{1}{\alpha +1} A &  -R^T\\
                R &  0
            \end{pmatrix}\right)
        \begin{pmatrix}
            \exp_m\left(t_j\left(  \frac{\alpha}{\alpha +1}A\right)\right)\\
            0
        \end{pmatrix} 
             $} %for inputs $t_jA, t_jR$ according to \eqref{eq:OrthMatpExp}
        }
 %\leftarrow \text{Exp4Geo}(t_jA, t_jR, \alpha)$}
   \hfill \COMMENT{cf. \eqref{eq:OrthMatpExp}}
    \ENDFOR\hfill \COMMENT{p-factor representation of geodesic $U^s(t_j) = UM(t_j) + QN(t_j)$}
    %\STATE{$\widetilde U^s \leftarrow\Exp_U(\Delta)$}  \hfill \COMMENT{hit of current shot}
    \STATE{ $A^s \leftarrow M(t_m)-\widehat M $, $R^s \leftarrow N(t_m)-\widehat M $}
            \hfill \COMMENT{p-factors of current gap vector}
    \STATE{$\gamma \leftarrow \sqrt{\|A^s\|^2 + \|R^s\|^2}$}
    %
    %%
    %%% parallel transport
    \FOR{$j=m,\ldots,0$}
    \STATE{ 
%        $\begin{pmatrix}
%               A^s\\ R^s 
%              \end{pmatrix}
     $[A^s,  R^s]    = $ ParaTrans\_pFactors$(M(t_j), N(t_j), A^s, R^s, \gamma)$
        }
    \hfill \COMMENT{initial projection plus approximate parallel transport}
    \ENDFOR
    \STATE{ update $A\leftarrow A - A^s, \quad R\leftarrow R - R^s$}
    \hfill \COMMENT{updated $\Delta=UA+QR$}
  \ENDWHILE
  \ENSURE{$\Delta = UA + QR$}
\end{algorithmic}
\end{algorithm}
\tc{A few remarks on this procedure are in order:
Observe that step 8} of \Cref{alg:p-Shooting} is the only step that depends on the chosen metric and thus makes the only difference when computing the Stiefel logarithm for the canonical, the Euclidean or any other $\alpha$-metric.
The subroutine ParaTrans\_pFactors$(M(t), N(t), A^s, R^s, \gamma)$ that appears in step 13 of \Cref{alg:p-Shooting}
realizes an approximation of the parallel transport of 
$\Delta^s = UA^s + QR^s$ along the geodesic using the unique representation with the $p$-factors $A^s, R^s$. This subroutine is detailed in Algorithm \ref{alg:ParaTrans_p}.
\begin{algorithm}
\caption{ParaTrans\_pFactors: Map the tangent vector $\Delta=UA_1 + QR_1\in T_{U_1}$ to
the tangent space $T_{U_2}$, preserve the length}
\label{alg:ParaTrans_p}
\begin{algorithmic}[1]
  \REQUIRE{$M_2,N_2,A_1,R_1\in \R^{p\times p}$, $\gamma > 0$\\
  \COMMENT{\footnotesize here, $M_2,N_2$ represent $U_2=UM_2+QN_2\in St(n,p)$,
  $A_1,R_1$ represent $\Delta = UA_1+QR_1\in T_{U_1}$.}}
  \STATE{ $S \leftarrow \sym(M_2^TA_1 + N_2^TR_1)$}
  \STATE{ $A_2 = A_1 - M_2S$, $R_2 = R_1 - N_2S$}
  \STATE{ $l =\sqrt{\|A_2\|^2 + \|R_2\|^2}$}
  \IF{ $l >\epsilon$}
    \STATE{$A_2 = \frac{\gamma}{l}A_2$, $R_2 = \frac{\gamma}{l}R_2$}\hfill \COMMENT{\footnotesize rescale to original length}
  \ELSE
    \STATE{$A_2 =0$, $R_2=0$.}
  \ENDIF
  \ENSURE{$A_2, R_2$}

This algorithm executes the same operation as in step 10 of \Cref{alg:SimpleShooting}
but on the representative $p\times p$ matrix factors.
\end{algorithmic}
\end{algorithm}
A numerical experiment that illustrates the performance of \Cref{alg:p-Shooting}
is given in Section \ref{sec:NumEx}.
%
%%
%%%
%%%%
%%%%%
%%%%%%
%%%%%%%
\subsection{An improved algebraic Stiefel logarithm for the canonical metric}
\label{sec:AlgImpro}
Let $U,\widetilde U=UM+QN\in St(n,p)$ with $M= U^T\widetilde U, QN = (I-UU^T)\widetilde U$ 
and an orthonormal completion
$\begin{pmatrix}M & X_0\\ N & Y_0\end{pmatrix}\in SO(2p)$
as introduced in the previous sections.
In \cite{StiefelLog_Zimmermann2017}, it has been shown that for solving
\eqref{eq:OrthMatpLog1} in the case of the canonical metric, it is sufficient to 
find $\Gamma\in \Skew(p)$ such that
\begin{equation}
 \label{eq:fundamentalMatrixEqCanLog}
 \begin{pmatrix}[c|c] 0 & I_p\end{pmatrix}
     \log_m\left(
              \begin{pmatrix}M & X_0\\ N & Y_0\end{pmatrix}
              \begin{pmatrix}I_p & 0\\ 0 & \exp_m(\Gamma)\end{pmatrix}
            \right)
     \begin{pmatrix}0 \\ I_p\end{pmatrix}
     = 0 \in \R^{p\times p}. 
\end{equation}
With 
$V=\begin{pmatrix}M & X_0\\ N & Y_0\end{pmatrix} = \exp_m \begin{pmatrix}A_0 & -B^T_0\\ B_0 & C_0\end{pmatrix}$
and $W=\exp_m\begin{pmatrix} 0 & 0\\ 0 & \Gamma\end{pmatrix}$,
\eqref{eq:fundamentalMatrixEqCanLog} requires to find $\Gamma$ such that the lower $p$-by-$p$ diagonal block
of the matrix $\log_m(VW)$ vanishes.
\tc{Let $[V,W] = VW-WV$ denote the matrix commutator.} The Baker-Campbell-Hausdorff (BCH) series
for the matrix logarithm is
    \begin{align}
    \nonumber
 &\log_m(VW)
       = \log_m(V) + \log_m(W)+ \frac{1}{2}[\log_m(V), \log_m(W)] \\
    \nonumber
   &  +\frac{1}{12}
      \left(
        \bigl[\log_m(V), [\log_m(V),\log_m(W)]\bigr]
       + \bigl[\log_m(W), [\log_m(W), \log_m(V)]\bigr]
      \right)+ \ldots,
    \end{align}
see \cite[\S 1.3, p. 22]{rossmann2006lie}.
The algorithm \cite[Alg. 1]{StiefelLog_Zimmermann2017} and the associated convergence analysis 
rely on the fact that with the choice of $\Gamma_0 = -C_0$,
the lower $p$-by-$p$ block of $\log_m(VW)$ vanishes up to terms of third order in the BCH series.
The algorithm iterates on this observation and produces the  sequence
\begin{equation}
 \label{eq:basicStLog}
 \begin{pmatrix}A_{k+1} & -B^T_{k+1}\\ B_{k+1} & C_{k+1}\end{pmatrix}
 := \log_m\left(
 \exp_m\begin{pmatrix}A_{k} & -B^T_{k}\\ B_{k} & C_{k}\end{pmatrix}
 \exp_m\begin{pmatrix}0 & 0\\ 0 & \Gamma_k\end{pmatrix}\right),
\end{equation}
with $\Gamma_k = -C_k$. It is guaranteed that $\|C_k\|\to 0$ for $k\to\infty$ at a linear rate as long as the input points $U,\widetilde U$ are close enough\tc{, see \cite[Theorem 4.1]{StiefelLog_Zimmermann2017}.}

In fact, up to commutator products of order three or higher,
the BCH series expansion of the lower diagonal block of
$\log_m\left( V_k W\right)
$
is
\begin{align*}
 C_k &+ \Gamma + \frac12(C_k \Gamma - \Gamma C_k)\\
 +& \frac{1}{12}
   \left(
     C_k^2\Gamma + \Gamma C_k^2 +C_k\Gamma^2 + \Gamma^2C_k-(B_kB_k^T\Gamma + \Gamma B_kB_k^T)
     -2(C_k\Gamma C_k + \Gamma C_k\Gamma)
   \right)\\
 =&C_k + \Gamma + \frac12(C_k \Gamma - \Gamma C_k)
  - \frac{1}{12}
   \left(
     B_kB_k^T\Gamma + \Gamma B_kB_k^T
   \right) + \text{h.o.t.}
\end{align*}
where ``$\text{h.o.t.}$'' comprises the higher-order terms in the BCH series.
Ignoring the quadratic terms $\Gamma C_k, C_k\Gamma$ and the higher ones and \tc{setting the above expression to zero} yields
a symmetric Sylvester equation for $\Gamma$:
\begin{equation}
\label{eq:sylv4log}
  C_k =S_k \Gamma + \Gamma S_k, \quad \text{ with } S_k:= \left(\frac{1}{12}B_kB_k^T-\frac12I_p\right).
\end{equation}
A sufficient criterion that guarantees a unique solution is that
$\| B_k\|_2 < \sqrt{6}$, since in this case
$\frac{1}{6}\|B_kB_k^T\|_2<1$, which entails that all eigenvalues of
$S_k$ are strictly negative. This in turn yields that 
$S_k$ and $-S_k$ have disjoint spectra, which ensures the unique solvability of \eqref{eq:sylv4log}, \cite[Section VII.2]{Bhatia1997}.
When selecting $\Gamma_k$ as the solution to \eqref{eq:sylv4log} at each iteration $k$, then the lower $p$-by-$p$ diagonal block of the BCH series vanishes up to fourth order commutator terms 
and terms that are quadratic in $C_k$ and $\Gamma_k$.
This idea was first presented in \cite{Zimmermann_Oberwolfach2018}.
\begin{remark}
This choice does not cancel all terms that are of first order in $\Gamma$ in the BCH series of
 $\log_m\left(
 V
  W\right)
$.
The series is ordered by the degree of nested commutator brackets.
In the $j$th-order commutator term, ``words'' formed by  $j$ ``letters'' of the two-letter alphabet 
$\{ \log_m(V), \log_m(W)\}$ appear, where each of $\log_m(V), \log_m(W)$ appears at least once, see 
\cite{Thompson1989, Thompson1989b, VanBruntVisser2015}.
This means that no matter where we cut off the tail of the BCH series, there remain terms that are linear in $\|\Gamma\|$ (or $\|C_k\|$ for that matter).
Expanding this series for $\log_m\left(
 VW\right)
$ in powers of $W$ as in \cite[Section 8.1]{Mueger2019} will not solve the issue,
because the collection of all terms that are of first order in $\|W\|$ constitutes an infinite series by itself that needs to be truncated in numerical applications.
\end{remark}
%
%
%
%
% The experiments presented in Section \ref{sec:NumEx} show that choosing $\Gamma_k$ 
% according to \eqref{eq:sylv4log} at every iteration step of \cite[Alg. 1]{StiefelLog_Zimmermann2017}
% improves the iteration count  by a factor of $2$ and accelerates the algorithm by a factor of $1.5$.
%
In the numerical implementation of \cref{alg:StLogBCH}, we introduce a Boolean ``Flag\_Sylv on/off'' to switch between the original version \cite[Alg. 1]{StiefelLog_Zimmermann2017}, which works with
$\Gamma_k = -C_k$  and the ``Sylvester-enhancement'', which selects $\Gamma_k$ as the solution to \eqref{eq:sylv4log}.
\begin{algorithm}
\caption{Improved algebraic Stiefel logarithm, canonical metric ($\alpha=0$).}
\label{alg:StLogBCH}
\begin{algorithmic}[1]
  \REQUIRE{$U, \widetilde{U}\in St(n,p)$, $\epsilon>0$
  convergence threshold, Boolean ``Flag\_Sylv on/off''}
  \STATE{ $M:= U^T\widetilde{U} \in \R^{p\times p}$}
  \STATE{ $QN := \widetilde{U} - UM \in \R^{n \times p}$}
    \hfill \COMMENT{compact QR-decomp.}
  \STATE{ $V_0 := \begin{pmatrix}M&X_0\\N&Y_0\end{pmatrix} \in O_{2p\times 2p}$}
    \hfill \COMMENT{orthonormal completion}

  \FOR{ $k=0,1,2,\ldots$}
  \STATE{$\begin{pmatrix}A_k & -B_k^T\\ B_k & C_k\end{pmatrix} := \log_m(V_k)$}
    \hfill \COMMENT{matrix log, $A_k, C_k$ skew}
  \IF{$\|C_k\|_2 \leq \epsilon$}
    \STATE{break}
  \ENDIF
  
  \IF{Flag\_Sylv}
    \STATE{$S_k := \frac{1}{12}B_kB_k^T -\frac12 I_p$}
    \STATE{solve $C_k = S_k \Gamma + \Gamma S_k$ for $\Gamma$}
       \hfill \COMMENT{sym. Sylvester equation}
  \ELSE
    \STATE{$\Gamma := -C_k$}\hfill \COMMENT{cancel first term in BCH series}
  \ENDIF
  
  \STATE{ $\Phi_{k} := \exp_m{(\Gamma)}$}
    \hfill \COMMENT{matrix exp, $\Phi_{k}$ orthogonal}
  \STATE{$V_{k+1} := V_{k}W_k$, where $W_k:=\begin{pmatrix}I_p& 0\\ 0 & \Phi_k\end{pmatrix}$
  }
    \hfill \COMMENT{update}
  \ENDFOR
  \ENSURE{$\Delta := \Log_{U}^{St}(\widetilde{U}) = U A_k + QB_k \in T_{U}St(n,p)$}
%
%  \textbf{Note:} Rentmeesters arrives at essentially the same 
%  algorithm \cite[Alg. 4, p.~91]{Rentmeesters2013} but from the different 
%  perspective of Riemannian optimization.
\end{algorithmic}
\end{algorithm}
A detailed convergence analysis can be conducted as in \cite{StiefelLog_Zimmermann2017} but is not worthwhile in the context of this work. Yet, the next proposition allows to compare the asymptotic convergence rate of \Cref{alg:StLogBCH} with ``Flag\_Sylv'' switched off, which is \cite[Alg. 1]{StiefelLog_Zimmermann2017} 
and \Cref{alg:StLogBCH} with ``Flag\_Sylv'' switched on, which is based on solving  \eqref{eq:sylv4log}.
\begin{proposition}
 \label{prop:BCHlog}
 Suppose that the input data $U\neq\widetilde{U}\in St(n,p)$ are such that \Cref{alg:StLogBCH} converges with ``Flag\_Sylv'' switched on.
 Assume further that there is a bound $0<\delta<1$ such that for $\log_m(V_{k}) =  \begin{pmatrix}A_k & -B_k^T\\ B_k & C_k\end{pmatrix}$, it holds $\|\log_m(V_k)\|_2< \delta$ throughout the algorithm's iteration loop.\footnote{For the original Stiefel log algorithm,
 conditions for the existence of such a global bound $\delta$ are established in \cite[Lemma 4.4]{StiefelLog_Zimmermann2017}.
 Similar techniques apply in the present context.
 }
 Then, for $k$ large enough, it holds 
 \[
  \|C_{k+1}\|_2\leq \left(\frac{6}{6-\delta^2} \frac{\delta^4}{1-\delta} + \mathcal{O}(\|C_{k}\|_2)\right) \|C_{k}\|_2.
 \]
 This implies the asymptotic convergence rate of \Cref{alg:StLogBCH} with ``Flag\_Sylv'' switched on for $k\to \infty$.  
 This compares to the asymptotic rate of \cite[Alg. 1]{StiefelLog_Zimmermann2017}, which according to 
 \cite[eq. (12)]{StiefelLog_Zimmermann2017} is bounded by
  \[
  \|C_{k+1}\|_2\leq \left(\frac{1}{6} \delta^2 + \frac{\delta^4}{1-\delta} + \mathcal{O}(\|C_{k}\|_2)\right) \|C_{k}\|_2.
 \]
\end{proposition}
\begin{proof}
 Since $\|B_k\|_2 \leq \|\log_m(V_k)\|<\delta < 1$,
 the matrices $S_k = \frac{1}{12}B_kB_k^T -\frac12 I_p$ are negative definite with largest eigenvalue bounded by $-\frac12 +\frac{\delta^2}{12} = -\frac{6-\delta^2}{12}$. Hence, the spectra of the symmetric matrices $S_k$ and $-S_k$ are separated by a vertical strip of width $\frac{6-\delta^2}{6}$ in the complex plain.
 Applying \cite[Theorem VII.2.12]{Bhatia1997} to the Sylvester equation
 $S_k\Gamma + \Gamma (-S_k) = C_k$ yields
 $\|\Gamma\|_2 \leq \frac{6}{6-\delta^2} \|C_k\|_2$.
 Calculating $C_{k+1}$ according to \eqref{eq:basicStLog} but with this choice of $\Gamma$ shows that all terms up to order four in the BCH series are at least quadratic in $\|C_k\|_2$:
 \begin{eqnarray*}
  C_{k+1} &=& \frac12 [C_k,\Gamma] + \frac{1}{12}([C_k^2,\Gamma] + [C_k,\Gamma^2]) - 2(C_k\Gamma C_k + \Gamma C_k\Gamma)\\
  && +\frac{1}{24} \left([B_kB_k^T,\Gamma^2]- (C_k^2\Gamma^2 + \Gamma^2C_k^2) +2[C_k\Gamma C_k,\Gamma]\right) + \text{h.o.t.}(5),
 \end{eqnarray*}
where $\text{h.o.t.}(5)$ are the terms of fifth order and higher
in the BCH series.
From $U\neq \widetilde{U}$, we get $\lim_{k\to \infty}\begin{pmatrix}A_k & -B_k^T\\ B_k & C_k\end{pmatrix} \neq 0$.  
Hence, for $k$ large enough, it holds $\|\Gamma\|_2 \leq  \frac{6}{6-\delta^2} \|C_k\|_2 \leq \|\begin{pmatrix}A_k & -B_k^T\\ B_k & C_k\end{pmatrix}\|_2 =  \|\log_m(V_k)\|_2$.
Because $C_{k+1}$ is but the lower diagonal subblock of 
 $\log_m\left(V_k W\right)$, the higher-order terms are bounded by
$\|\text{h.o.t}(5)\|_2 \leq \sum_{l=5}^\infty \|\log_m(V_k)\|_2^{l-1} \|\Gamma\|_2\leq \frac{6}{6-\delta^2} \|C_k\|_2 \frac{\delta^4}{1-\delta}$, see \cite[Lemma A.1]{StiefelLog_Zimmermann2017}.
The claim is now a straightforward consequence.
\end{proof}
\Cref{prop:BCHlog} shows that with smaller values of $\delta$, the Sylvester approach becomes more and more favorable.
The actual value of $\delta$ depends on how close the inputs $U,\widetilde U$ are.
For $\delta \approx 0.7147$ %0.714740520345120
the bounds for the asymptotic convergence rates are the same for both the  approaches ``Flag\_Sylv on/off''; for $\delta \approx0.3286$ %0.328559198972538
the rate of the Sylvester-based method is improved a factor of $2$
and by a factor of $10$ for $\delta \approx 0.1270.$ %0.126961175098983
Note that in both cases, the bounds overestimate the true convergence rates.

The main computational effort of \Cref{alg:StLogBCH} is in the computation of the matrix logarithm in step 5.
This can be achieved efficently (and without resorting to complex numbers arithmetics) by first computing a real Schur form. Assuming that the principal matrix logarithm is properly defined, the Schur form of an orthogonal matrix is block-diagonal with blocks $(1)$ of size $(1\times 1)$ and 
$(2\times 2)$-blocks of the form
 $\begin{pmatrix}
  \cos(\varphi) & -\sin(\varphi)\\
  \sin(\varphi) & \cos(\varphi)
  \end{pmatrix}
$, the matrix logarithm of such a block being 
 $\begin{pmatrix}
  0 & -\varphi\\
 \varphi & 0
  \end{pmatrix}
$. This is exploited in the actual implementation. 
%
%%
%%%
%%%%
%%%%%
%%%%%%
%%%%%%%
\subsection{A geodesic Newton method for the Stiefel logarithm}
\label{sec:GeoNewton}
The nonlinear matrix equation \eqref{eq:OrthMatpLog1} can be cast in the following form
\[
 V^T
 F \begin{pmatrix}
  A & -B^T\\
  B & C
 \end{pmatrix}
 = I, \quad 
 \text{where }
 V =   \begin{pmatrix}
  M & X_0\\
  N & Y_0
 \end{pmatrix}\in SO(2p).
\]
For brevity, introduce 
\[
 \widehat F: \Skew(2p) \to SO(2p), \hspace{0.2cm} S\mapsto V^TF(S).
 \] 
%  \quad   \text{where }
%  S = 
%  \begin{pmatrix}
%   A & -B^T\\
%   B & C
%  \end{pmatrix},.
% \]
The parameter domain $\Skew(2p)$ has a vector space structure that allows to employ Euclidean techniques. The co-domain, however is the Lie group $SO(2p)$.
Ignoring the structure of $SO(2p)$, the classical Newton method can be applied 
to the root finding problem $V^TF(S) -I = 0$.
Given a starting point $S_0$, the Newton method requires to solve the Taylor-linearized problem
\[
 I = \widehat F(S)+ D\widehat F_S(H) \approx \widehat F(S +H) 
\]
so that the update $H$ is determined by the linear system
$D\widehat F_S(H) = I-\widehat F(S)$.
Yet, the left-hand side and the right-hand side are not compatible, as 
$\widehat F(S) \in SO(2p)$ and 
$D\widehat F_S(H) \in T_{\widehat F(S)}SO(2p) = \widehat F(S) \Skew(2p)$.

A remedy is to replace the Euclidean first-order Taylor approximation,
which can be thought of as progressing along a straight line, by moving along a geodesic in the same direction. 
% (whose first-order Taylor approximation coincides with the Euclidean Taylor). 
This leads to a first-order approximation that preserves the Riemannian structure
\[
  \widehat F(S +H)\approx  \Exp_{\widehat F(S)}(D\widehat F_S(H)).
\]
On $SO(2p)$, the geodesic that starts from $\widehat F(S)$ in 
the direction of $D\widehat F_S(H)$ is
\[
 \Exp_{\widehat F(S)}(D\widehat F_S(H)) = \widehat F(S) \exp_m(\widehat F(S)^T D\widehat F_S(H)),
\]
see, e.g., \cite{godement2017introduction}.
The first-order approximation to  \eqref{eq:OrthMatpLog1} becomes
\begin{equation}
 \label{eq:GeoNewt}
 I = \widehat F(S) \exp_m(\widehat F(S)^T D\widehat F_S(H)) \quad
 \Leftrightarrow\quad
  \underbrace{\widehat F(S)^T D\widehat F_S(H)}_{\text{ skew}} = \underbrace{\log_m(\widehat F(S)^T)}_{\text{ skew}}.
\end{equation}
The left-hand side is a linear operator
\[
 L_S: \Skew(2p)\to \Skew(2p), \hspace{0.2cm} H\mapsto \widehat F(S)^T D\widehat F_S(H).
\]
In summary, this leads to the iterative scheme stated in \Cref{alg:GeoNewt}.
Upon convergence,
$S_k = 
 \begin{pmatrix}
  A_k & -B^T_k\\
  B_k & C_k
 \end{pmatrix}
 $
 is found and the output
 $\Delta := \Log_{U}^{St}(\widetilde{U}) = U A_k + QB_k \in T_{U}$
 is formed with $U,Q$ as in \Cref{alg:StLogBCH}.
Computationally, the expensive part is to evaluate the operator $L_{S}$.
For convenience, let us restrict to the Euclidean metric and
write 
\[
 S = \begin{pmatrix}
  A & -B^T\\
  B & C
 \end{pmatrix},\quad
 S_{\text{tri}} = 
 \begin{pmatrix}
  2A & -B^T\\
  B & 0
 \end{pmatrix},
 \quad
 S_{\text{di}}
 =\begin{pmatrix}
  A & 0\\
  0 & C
 \end{pmatrix},
\]
likewise for $H\in \Skew(2p)$.
It holds
\begin{eqnarray}
\nonumber
 L_S(H) &=& \widehat F(S)^T D\widehat F_S(H) = F(S)^TV V^T DF_S(H)\\
 \label{eq:GeoNewtOperator}
 &=& F(S)^T \left(D(\exp_m)_{S_{\text{tri}}}(H_{\text{tri}})\exp_m(-S_{\text{di}})
 + \exp_m(S_{\text{tri}})D(\exp_m)_{S_{\text{di}}}(-H_{\text{di}})\right).
\end{eqnarray}
In an implementation, we do not actually form this operator, but implement its action on a matrix $H$. This is then used in a matrix-free version of the GMRES algorithm \cite{GMRES_Saad_1986}
(as pre-installed in Matlab, version R2019b).
% Formally, the derivative of the matrix exponential is
% \[
%  D(\exp_m)_S(H) = \exp_m(S)(H -\frac1{2!}[S,H] + \frac{1}{3!} [S,[S,H]]\pm\ldots),
% \]
% see \cite[\S 5.4]{Hall_Lie2015}.
For the practical evaluation of  $D(\exp_m)_S(H)$, we use Mathias' Theorem \cite[Thm 3.6]{Higham:2008:FM}, which simultaneously gives $\exp_m(S)$.
\begin{algorithm}
 \caption{GeoNewton for solving \eqref{eq:OrthMatpLog1}}
\begin{algorithmic}[1]
 \STATE{$k=0$}
 \WHILE{ $\| \log_m(\widehat F(S_k)^T)\|_{F}>\tau$} 
 \STATE{solve  $L_{S_k}(H_k) = \log_m(\widehat F(S_k)^T)$}
 \STATE{update $S_{k+1} = S_k + H_k$}
 \STATE{$k=k+1$}
 \ENDWHILE 
\end{algorithmic}
\label{alg:GeoNewt}
\end{algorithm}
An equivalent alternative to \eqref{eq:OrthMatpLog1} is to solve
\begin{equation}
\label{eq:EucNewton4Log}
 0 = \log_m(F(S)) - \log_m(V).
\end{equation}
In this form, both the unknown $S$ and the output $\log_m(F(S)) - \log_m(V)$
are in the vector space $\Skew(2p)$ so that this equation is amenable to be treated with the classical Newton method.
However, computing or approximating the derivative remains the computational bottleneck. The associated linear operator now involves differentials of both the matrix exponential and the matrix logarithm. We tackle this with the same strategy as above by relying on Mathias' Theorem and the matrix-free GMRES method.
The matrix functions $\exp_m(S)$ and  $\log_m(S)$ may be approximated via the Cayley transformation $\Cay(S) = (I-\frac12 S)^{-1} (I + \frac12 S)$ and its inverse;
likewise the differentials of $D(\exp_m)_S(H)$ and  $D(\log_m)_S(H)$ may be approximated with the differentials of the corresponding Cayley transformations.

We mention these approaches for the sake of completeness
and include the algorithms based on \eqref{eq:GeoNewt} and \eqref{eq:EucNewton4Log}  in the algorithmic competition.
However, the numerical experiments show that none of the above approaches can't compete, neither with the $p$-shooting method 
\Cref{alg:p-Shooting} nor with the algebraic Stiefel log algorithm
\Cref{alg:StLogBCH}. This also holds, when the Cayley transformations are used to replace $\exp_m$ and $\log_m$. Therefore, we omit a detailed discussion.
%
%
%
%
%
%
%
%
%%
%%%
%%%%
%%%%%
%%%%%%
%%%%%%%
\section{Numerical experiments}
\label{sec:NumEx}
In this section, we conduct various numerical experiments for assessing the performance of the proposed approaches to solve the local geodesic endpoint problem.
All experiments are performed with Matlab R2019b on a Linux 64bit HP notebook with 
four Intel(R) Core(TM) i7-5600U 2.60GHz CPUs.\footnote{The Matlab code and an accompanying Python implementation is available on \\
\begin{scriptsize}
\url{https://github.com/RalfZimmermannSDU/RiemannStiefelLog/tree/main/Stiefel_log_general_metric/}
\end{scriptsize}
}
% \tc{An additional numerical test with more practical relevance is included \Cref{sec:TriExperiment} of the supplements, where we employ the Riemannian log-methods to compute the Riemannian center of mass of a geodesic triangle on $St(n,p)$.}

\subsection{The geodesic endpoint problem for the canonical metric}
\label{sec:1Experiment}
%
% TABLE DATA
%
%
%n=2000$, $p=500$, tau = 1.0e-11, dist = 5pi, 5 runs
% The average timing of the various methods is:
% 
% ans =
% 
%    13.0036
%     8.2914
%     8.2773
%    25.5299
%    53.8908
%    34.0636
% 
% The average iteration count of the various methods is:
% 
% ans =
% 
%    13.0000
%     7.0000
%     7.0000
%    41.0000
%    35.0000
%    39.8000
% 
% The average reconstruction accuracy of the various methods is:
% 
% ans =
% 
%    1.0e-11 *
% 
%     0.5001
%     0.0288
%     0.0293
%     0.4862
%     0.6132
%     0.5062
%%
%%%
%%%%
%%%%%
%%%%%%******table***********table************table***********************table*******************
\begin{table}[ht]
\begin{small}
\begin{center}
\begin{tabular}{l|l|l|l}
 \rowcolor{gray!20}
 \hline
  \multicolumn{4}{c}{(Sec. \ref{sec:1Experiment}) {\bf Test Case 1:} random data $n=2000$, $p=500$, canonical metric, 5 runs}\\ 
  \hline
  \multicolumn{4}{c}{ $\dist(U,\widetilde U) = 5\pi$ }\\%\smallskip\\
  \hline  \rowcolor{gray!20}
  Method      & av. rel. error $\|\Delta - \Delta_{rec}\|_{\infty}$ & av. iter. count & av. time  \\

 \hline 
  Alg. \ref{alg:StLogBCH}                  & $0.50\cdot 10^{-11}$  & $13.0$ &  $13.00$s\\
 \hline
 \rowcolor{gray!10}
 Alg. \ref{alg:StLogBCH}+Sylv.           & $0.29\cdot 10^{-12}$ & $7.0$  & $8.29$s\\
 \hline
 Alg. \ref{alg:StLogBCH}+Sylv.+Cay.    & $0.29\cdot 10^{-12}$ & $7.0$  & $8.27$s\\
  \hline
 \rowcolor{gray!10}
 Alg. \ref{alg:p-Shooting} on 2 steps      & $0.49\cdot 10^{-11}$  & $41.0$ & $25.53$s\\
 \hline
 Alg. \ref{alg:p-Shooting} on 4 steps      & $0.61\cdot 10^{-11}$  & $35.0$ & $53.89$s\\
  \hline
    \rowcolor{gray!10}
  Alg. \ref{alg:SimpleShooting} on 2 steps & $0.51\cdot 10^{-11}$  & $39.8$ & $34.06$s\\
  \hline
     Single shooting \cite[\S 2.3]{Sutti_PhD_2020}  &   
     \multicolumn{3}{c}{\scriptsize unfeasible du to memory overflow}\\
 \hline 
  \multicolumn{4}{c}{ }\\%\smallskip\\
%%
%%
%%
%%
%%$n=120$, $p=30$, tau = 1.0e-11, dist = pi
%%The average timing of the various methods is
% ans =
% 
%     0.0267
%     0.0175
%     0.0180
%     0.0336
%     0.0637
%     0.0322
% 
% The average iteration count of the various methods is:
% 
% ans =
% 
%    10.2000
%     5.0000
%     5.0000
%    26.8000
%    24.7000
%    26.4000
% 
% The average reconstruction accuracy of the various methods is:
% 
% ans =
% 
%    1.0e-11 *
% 
%     0.2259
%     0.1592
%     0.1604
%     0.2913
%     0.1934
%     0.2807
%%
%%%
 \hline
  \rowcolor{gray!20}
  \multicolumn{4}{c}{(Sec. \ref{sec:1Experiment}) {\bf Test Case 2:} random data $n=120$, $p=30$, canonical metric, 10 runs}\\
  \hline
  \multicolumn{4}{c}{ $\dist(U,\widetilde U) = \pi$ }\\%\smallskip\\
  \hline
  \rowcolor{gray!20}
  Method      & av. rel. error $\|\Delta - \Delta_{rec}\|_{\infty}$ & av. iter. count & av. time  \\
 \hline 
  Alg. \ref{alg:StLogBCH}                 & $0.226\cdot 10^{-11}$  & $10.2$ & $0.027$s\\
 \hline
 \rowcolor{gray!10}
 Alg. \ref{alg:StLogBCH}+Sylv.            & $0.159\cdot 10^{-11}$  & $5.0$  & $0.018$s\\
 \hline
 Alg. \ref{alg:StLogBCH}+Sylv.+Cay.       & $0.160\cdot 10^{-11}$  & $5.0$  & $0.018$s\\
  \hline
 \rowcolor{gray!10}
 Alg. \ref{alg:p-Shooting} on 2 steps     & $0.291\cdot 10^{-11}$  & $26.8$ & $0.033$s\\
 \hline
 Alg. \ref{alg:p-Shooting} on 4 steps     & $0.193\cdot 10^{-11}$  & $24.7$ & $0.064$s\\
 \hline
  \rowcolor{gray!10}
  Alg. \ref{alg:SimpleShooting} on 2 steps& $0.281\cdot 10^{-11}$  & $26.4$ & $0.032$s\\
  \hline
    Single shooting \cite[\S 2.3]{Sutti_PhD_2020}  &   $0.58\cdot 10^{-14}$ ${ }_{\text{(results for 1 run.)}}$  & $5$& $524.6$s\\
 \hline
  \multicolumn{4}{c}{ }\\%\smallskip\\
%%
%%
%%
% Sutti, single shooting on St(120,30), dist=1.0 pi, tau = 1.0e-11
% 1 runs. 
%
% Timing 524.6279s = %
% Iters count: 5
% accuracy:    5.8441e-15
%%
%%
%%
%%
%%
%%

%%
%%%
 \hline
  \rowcolor{gray!20}
  \multicolumn{4}{c}{(Sec. \ref{sec:1Experiment}) {\bf Test Case 3:} random data $n=12$, $p=3$, canonical metric, 100 runs}\\
  \hline
  \multicolumn{4}{c}{ $\dist(U,\widetilde U) = 0.95\pi$ (averaging only over the converged runs) }\\%\smallskip\\
  \hline
  \rowcolor{gray!20}
  Method      & av. rel. error $\|\Delta - \Delta_{rec}\|_{\infty}$ & av. iter. count & av. time  \\
 \hline 
  Alg. \ref{alg:StLogBCH} &  $0.62\cdot 10^{-10}$ ${ }_{\text{(1 run not conv'd)}}$     & $120.3$ &  $0.046$s\\
 \hline
 \rowcolor{gray!10}
 Alg. \ref{alg:StLogBCH}+Sylv. & $0.50\cdot 10^{-10}$ ${ }_{\text{(1 run not conv'd)}}$ & $41.1$  & $0.023$s\\
 \hline
 Alg. \ref{alg:StLogBCH}+Sylv.+Cay.  & $0.53\cdot 10^{-10}$  ${ }_{\text{(1 run not conv'd)}}$& $41.3$ & $0.021$s\\
  \hline
 \rowcolor{gray!10}
 Alg. \ref{alg:p-Shooting} on 2 steps  & ${ }_{\text{(all 100 runs not conv'd)}}$  & -- & --\\
 \hline
 Alg. \ref{alg:p-Shooting} on 4 steps &   $0.80\cdot 10^{-10}$ ${ }_{\text{(all 100 runs conv'd.)}}$  & $212.2$& $0.031$s\\
 \hline
  \rowcolor{gray!10}
 Alg. \ref{alg:SimpleShooting} on 2 steps  & ${ }_{\text{(all 100 runs not conv'd)}}$  & -- & --\\
 \hline
  Single shooting \cite[\S 2.3]{Sutti_PhD_2020}  &   $0.42\cdot 10^{-14}$ ${ }_{\text{(41 runs not conv'd.)}}$  & $7.93$& $0.021$s\\
 \hline
%% Sutti, single shooting on St(12,3), dist=0.95 pi, tau = 1.0e-11
% converged runs. 59
%
% Timing
%  sum(time_array(1,entries))/59 
% ans =     0.0207
%
% Iters count: 7.9322
% accuracy:     4.2231e-15
%%
%%%
\end{tabular}
\caption{Numerical performance for the cases considered in \Cref{sec:1Experiment}.
}
\label{tab:numEx1}
\end{center}
\end{small}
\end{table}
%%%%%%***end******table*****************end******table********************************************
%%%%%
%%%%
%%%
%%
%
We create two points $U, \widetilde{U}$ pseudo-randomly on $St(n,p)$, but such that they are a prescribed geodesic distance apart. More precisely, we construct $U$ from a QR-decomposition of a random $(n\times p)$-matrix with entries sampled from the uniform distribution.
Then, we create a random tangent vector $\Delta = U A + (I-UU^T)T$, where $A\in \R^{p\times p}$ is random but skew and $T\in\R^{n\times p}$ is random. Then, $\Delta$ is scaled to the prescribed length according to the Riemannian metric and $\widetilde{U}\in St(n,p)$ is obtained as $\widetilde U = \Exp_U(\Delta)$.

Then, the various Stiefel-log algorithms are applied to compute the reconstructed tangent vector $\Delta_{\text{rec}}$ and the absolute accuracy is checked in the infinity-matrix norm. In summary, we perform the following three steps
\[
 %\frac{\|\Delta - \Delta_{rec}\|_{\text{fro}}}{\|\Delta\|_{\text{fro}}}.
(a)\quad  \widetilde U \leftarrow \Exp_U(\Delta), \quad (b)\quad \Delta_{\text{rec}} \leftarrow \Log_U(\widetilde U), \quad (c)\quad \text{compute } \|\Delta - \Delta_{\text{rec}}\|_{\infty}.
\]
We perform 10 runs with random data, record the numerical accuracy, the iteration count and the computation time and average over the results.
In order to evaluate the logarithm, we use the following methods:
\begin{itemize}
\item Alg. \ref{alg:StLogBCH} as in \cite{StiefelLog_Zimmermann2017}.
\item Alg. \ref{alg:StLogBCH} enhanced by the Sylvester equation \eqref{eq:sylv4log} as detailed in \Cref{sec:AlgImpro}.
\item Alg. \ref{alg:StLogBCH} enhanced by \eqref{eq:sylv4log} and with the Cayley transformation replacing the matrix exponential in Step 15.
\item Alg. \ref{alg:p-Shooting} on two time steps $\{0.0, 1.0\}$ in the unit interval $[0.1]$.
\item Alg. \ref{alg:p-Shooting} on four time steps $\{0.0,0.\bar{3}, 0.\bar{6}, 1.0\}$ of the unit interval $[0,1]$.
\item \tc{Alg. \ref{alg:SimpleShooting} on two time steps $\{0.0, 1.0\}$ in the unit interval $[0.1]$.}
\item \tc{The single shooting method of \cite[Section 2.3]{Sutti_PhD_2020} that is based on a Newton root finding problem under the canonical metric.}
\end{itemize}
In each case, the convergence threshold is set to $\tau = 10^{-11}$.
\Cref{tab:numEx1} displays the results for data on $St(n=2000, p = 500)$ with $\dist(U,\widetilde U) = 5\pi$, for data on $St(n=120, p = 30)$ with $\dist(U,\widetilde U) = \pi$ and for data on $St(n=12, p = 3)$ with $\dist(U,\widetilde U) = 0.95\pi$.
Example plots of the convergence histories are shown in 
\Cref{supp:fig:St_canoni_2000_500_5pi} and \Cref{supp:fig:St_canoni_120_30_1pi} of the supplements, respectively.

The table shows that \Cref{alg:StLogBCH} with ``Flag\_Sylv on'' exhibits the best performance for the cases considered in regards of the computation time. 
In terms of the numerical accuracy it is outranked by the Newton-based single shooting method of \cite[Section 2.3]{Sutti_PhD_2020}, provided that the latter converges.
It can also be seen that a subdivision of the interval $[0,1]$ is not required for \Cref{alg:p-Shooting} in order to converge for the larger data sets under consideration. The dimensions and distance for the low-dimensional data set on $St(12,3)$ are chosen as in \cite{Sutti_V:2020b}, where the global geodesic endpoint problem is considered. 
Recall that the estimated injectivity radius of $St(n,p)$ under the canonical metric is at least $0.89\pi$ and most likely not larger, see \cite{Rentmeesters2013}.
In fact, the methods considered here are local by nature. \tc{As is to be expected, in the experiments, we observe convergence in some cases and divergence in other cases}, see \Cref{tab:numEx1} for details. The \Cref{alg:p-Shooting} on four time steps in $[0,1]$ converges in all cases considered, while it diverges in all cases, if only the boundary points of $[0,1]$ are considered in the discrete parallel transport.
\subsection{The geodesic endpoint problem for the Euclidean metric}
\label{sec:2Experiment}
In this section, we repeat the experiments of \Cref{sec:1Experiment} with exactly the same set-up, but for the Euclidean metric.
Since the algebraic Stiefel log algorithm \cref{alg:StLogBCH} is not available for metrics other than the canonical one, we juxtapose the following methods
\begin{itemize}
\item Alg. \ref{alg:GeoNewt}, ``GeoNewton'', the geodesic Newton method for \eqref{eq:OrthMatpLog1}.
\item Alg. ``EucNewton'', based on the classical Newton method for solving \eqref{eq:EucNewton4Log}.
\item Alg. \ref{alg:p-Shooting} on two time steps $\{0.0, 1.0\}$ of the unit interval $[0.1]$.
\item Alg. \ref{alg:p-Shooting} on four time steps $\{0.0,0.\bar{3}, 0.\bar{6}, 1.0\}$ of the unit interval $[0,1]$.
\end{itemize}
In each case, the convergence threshold is set to $\tau = 10^{-11}$.
\Cref{tab:numEx2} displays the results for data on $St(n=120, p = 30)$ with $\dist(U,\widetilde U) = \pi$
and for data on $St(n=2000, p = 500)$ with $\dist(U,\widetilde U) = 5\pi$.
In the latter case, only one random run is performed due to the large computation times for the  Newton methods.
Example plots of the convergence histories are shown in 
\Cref{fig:St_Euclid_2000_500_5pi} and \Cref{fig:St_Euclid_120_30_1pi} of the supplements, respectively.
%%
%%
%%
%%
%%
%% Newton test n=2000, p=500, d=5pi, tau = 1.0e-11
% The average timing of the various methods is:
% 
% ans =
% 
%    10.4607       p-shooting 2
%    14.7905       p-shooting 4
%   958.1609       GeoNewton
%   890.1257       EucNewton
% 
% The average iteration count of the various methods is:
% 
% ans =
% 
%     20
%     11
%     13
%      6
% 
% The average reconstruction accuracy of the various methods is:
% 
% ans =
% 
%    1.0e-11 *
% 
%     0.2635
%     0.3597
%     0.6190
%     0.2340
%%
%%
%%%
%%%%
%%%%%
%%%%%%******table***********table************table***********************table*******************
\begin{table}[!hb]
\begin{small}
\begin{center}
\begin{tabular}{l|l|l|l}
\hline
  \rowcolor{gray!20}
  \multicolumn{4}{c}{(Sec. \ref{sec:2Experiment}) {\bf Test Case 1:} random data $n=2000$, $p=500$, Euclidean metric, 1 run}\\
  \hline
  \multicolumn{4}{c}{ $\dist(U,\widetilde U) = 5\pi$ }\\%\smallskip\\
  \hline
  \rowcolor{gray!20}
  Method& rel. error $\|\Delta - \Delta_{rec}\|_{\infty}$ & iter. count & time  \\
 \hline 
  Alg. 5 GeoNewton                      & $0.61\cdot 10^{-11}$ & $13$    & $958.2$s \\
  \hline
 \rowcolor{gray!10}
 Alg. EucNewton                       & $0.23\cdot 10^{-11}$ & $6$  & $890.1$s\\
  \hline
 Alg. \ref{alg:p-Shooting} on 2 steps & $0.26\cdot 10^{-11}$ & $20$ & $10.5$s\\
 \hline
 \rowcolor{gray!10}
 Alg. \ref{alg:p-Shooting} on 4 steps & $0.36\cdot 10^{-11}$ & $11$  & $14.8$s\\
 \hline
  \multicolumn{4}{c}{ }\\%\smallskip\\
%%
%%%
 \hline
  \rowcolor{gray!20}
  \multicolumn{4}{c}{(Sec. \ref{sec:2Experiment}) {\bf Test Case 2:} random data $n=120$, $p=30$, Euclidean metric, 10 runs}\\
  \hline
  \multicolumn{4}{c}{ $\dist(U,\widetilde U) = \pi$ }\\%\smallskip\\
  \hline
  \rowcolor{gray!20}
  Method& av. rel. error $\|\Delta - \Delta_{rec}\|_{\infty}$ & av. iter. count & av. time  \\
 \hline 
  Alg. 5 GeoNewton                      & $0.11\cdot 10^{-11}$ & $8.0$    & $0.34$s\\
 \hline
 \rowcolor{gray!10}
 Alg. EucNewton                       & $0.071\cdot 10^{-11}$ & $5.0$  & $0.85$s\\
  \hline
 Alg. \ref{alg:p-Shooting} on 2 steps & $0.078\cdot 10^{-11}$ & $13.1$ & $0.016$s\\
 \hline
 \rowcolor{gray!10}
 Alg. \ref{alg:p-Shooting} on 4 steps & $0.12\cdot 10^{-11}$ & $9.0$  & $0.030$s\\
 \hline
%
%%
%%%
\end{tabular}
\caption{Numerical performance for the cases considered in \Cref{sec:2Experiment}.
}
\label{tab:numEx2}
\end{center}
\end{small}
\end{table}
%%%%%%***end******table*****************end******table********************************************
%%%%%
%%%%
%%%
%% St(120,30) Euclid, d=pi, tau = 1.0e-11
%The average timing of the various methods is:
% 
% ans =
% 
%     0.0157
%     0.0295
%     0.3367
%     0.8506
% 
% The average iteration count of the various methods is:
% 
% ans =
% 
%    13.1000
%     9.0000
%     8.0000
%     5.0000
% 
% The average reconstruction accuracy of the various methods is:
% 
% ans =
% 
%    1.0e-11 *
% 
%     0.0775
%     0.0121
%     0.1081
%     0.0714
%
%
%
%

%
%
%
The table shows that \Cref{alg:p-Shooting} without a subdivision of $[0,1]$
exhibits the best performance in terms of the computation time for the cases considered here. For the test case with data on 
$St(2000,500)$, the method is ca. $90$ times faster than the Newton-based methods and also much more memory efficient, since no large linear operators need to be constructed.
These computation times are representative for other values of $\alpha$ in the family of Riemannian metrics.
\subsection{Investigations on the parametric dependencies} 
\label{sec:3Experiment}
\tc{In this section, the performance of the proposed methods for computing
$\Log_U^\alpha(\widetilde U)$ with $U,\widetilde U \in St(n,p)$ is investigated under changes in the metric parameter $\alpha$, the Riemannian distance of the input points $\dist(U, \widetilde U)$ as well as the matrix dimensions $n$ and $p$ .
}

\begin{figure}[ht]
\centering
\includegraphics[width=1.0\textwidth]{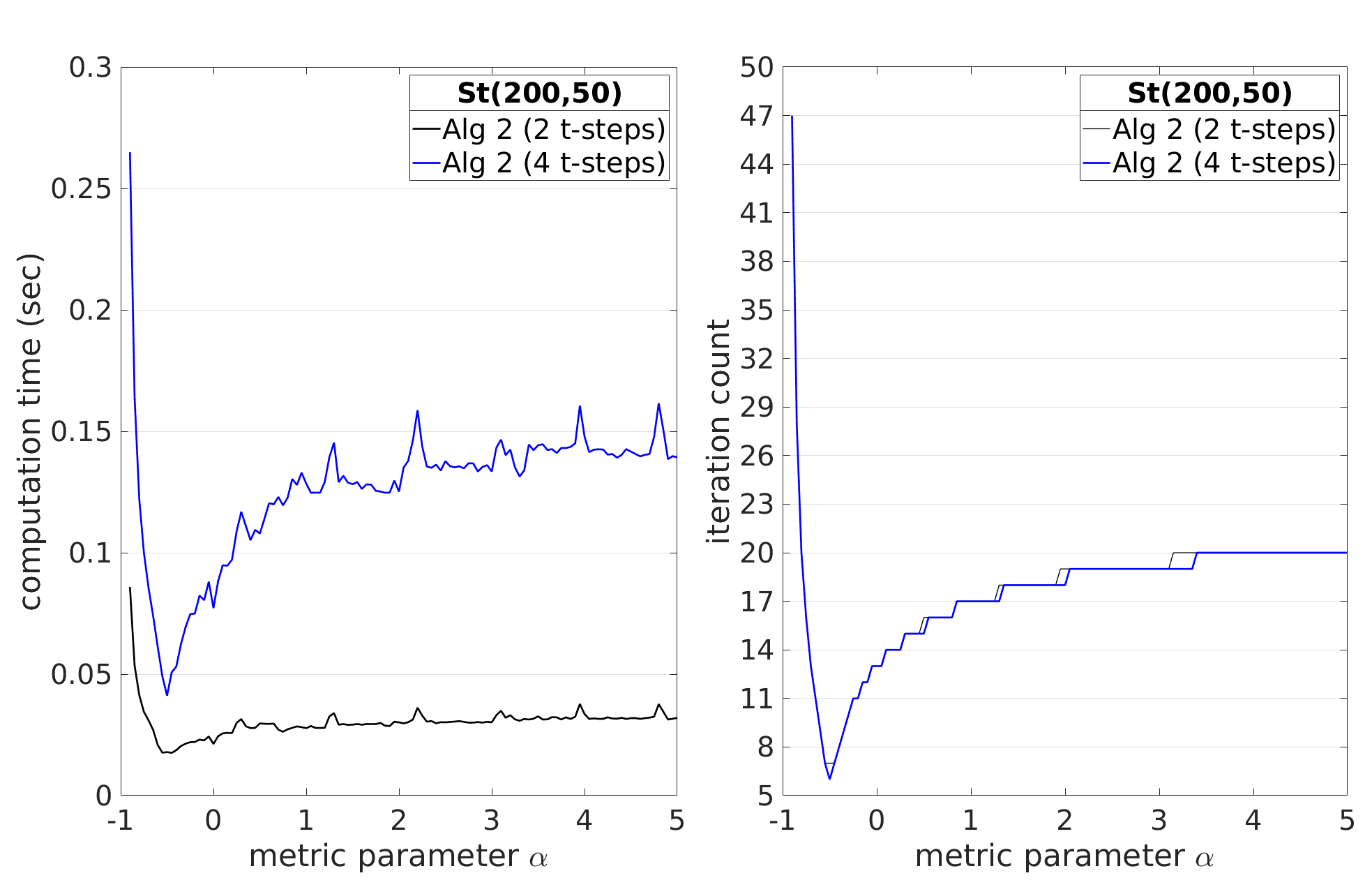}
\caption{(cf. \Cref{sec:3Experiment}) Wall clock computation time in seconds versus $\alpha$ (left) and iteration count until convergence versus $\alpha$ (right) for computing $\Log_U^\alpha(\widetilde U)$ with \Cref{alg:p-Shooting}. The trial parameter set is $\{\alpha=-0.9 + 0.05j|\hspace{0.1cm} j = 0,\ldots,118\}$.
The data points $U,\widetilde U\in St(200,50)$ are at a Riemannian $\alpha$-distance of $\dist_\alpha(U,\widetilde U) = 0.5\pi$. Timing results are averaged over 100 runs.
}
\label{fig:alpha_dependence_n200_p50_d_05pi_alpha_-09--5_av100}
\end{figure}

\tc{We start with assessing the performance of
\Cref{alg:p-Shooting} for metric parameters $\alpha\in (-1,\infty)$, where the associated metric is Riemannian.
(For $\alpha\in (-\infty, -1)$, the metric becomes pseudo-Riemannian, see \cite[Section 5.5]{HueperMarkinaLeite2020}.)
To this end, we fix the dimensions $n=200$, $p=50$ and construct pseudo-random data
$U\in St(200,50)$, $\Delta_0 \in T_USt(200,50)$ as described in \Cref{sec:1Experiment}.
We discretize the parameter interval $[-0.9, 5.0]$ with equidistant steps
of size $0.05$. For each value $\alpha\in \{-0.9 + 0.05j| j = 0,\ldots,118\}$, 
we compute $\Delta(\alpha) = \frac{d}{\|\Delta_0\|_\alpha} \Delta_0$, so that $\Delta(\alpha)$ is normalized to a length of $d$ according to the $\alpha$-metric. As a distance factor, we use $d=0.5\pi$.
Then, we set $\widetilde U = \Exp_U^\alpha(\Delta(\alpha))$ with the Stiefel exponential computed according to \eqref{eq:StAlphaRed}. In this way, a Stiefel data pair $U,\widetilde U$ with $\dist_\alpha(U,\widetilde U) = d = 0.5\pi$ is obtained. We apply \cref{alg:p-Shooting} to compute 
$\Log_U^\alpha(\widetilde U)$ up to a convergence threshold of $\tau = 10^{-11}$ and record the wall clock computation time as well as the iteration count.
The results are displayed in \Cref{fig:alpha_dependence_n200_p50_d_05pi_alpha_-09--5_av100}.
The supplements feature an analog experiment with data on $St(2000,200)$, see \Cref{fig:alpha_dependence_n2000_p200_d08pi_alpha_-09--5}.}
%
%

%
%---------------------------------------------------------------------------
\tc{In both test cases, the iteration count is minimal precisely at $\alpha = -\frac{1}{2}$, which corresponds to the Euclidean metric. We conjecture that the dependency of the iteration count on the metric parameter $\alpha$ is curvature-related. Recall that the curvature of a Riemannian manifold depends on the metric. Furthermore, observe that in the flat matrix space $\R^{n\times p}$, the shooting method produces the correct solution after one single iteration, because the geodesics are straight lines so that $\Exp_U(t\Delta) = U + t\Delta$.  Following this line of reasoning, the experimental results suggest that at least locally around the generic pseudo-random data points, $\left(St(n,p), \langle \cdot, \cdot\rangle^\alpha \right)$  is least curved for the induced Euclidean metric, where $\alpha = -\frac12$.
}
%
%
%%
%%%
%%%%
%%%%%
%%%%%%

\begin{figure}[ht]
\centering
\includegraphics[width=1.0\textwidth]{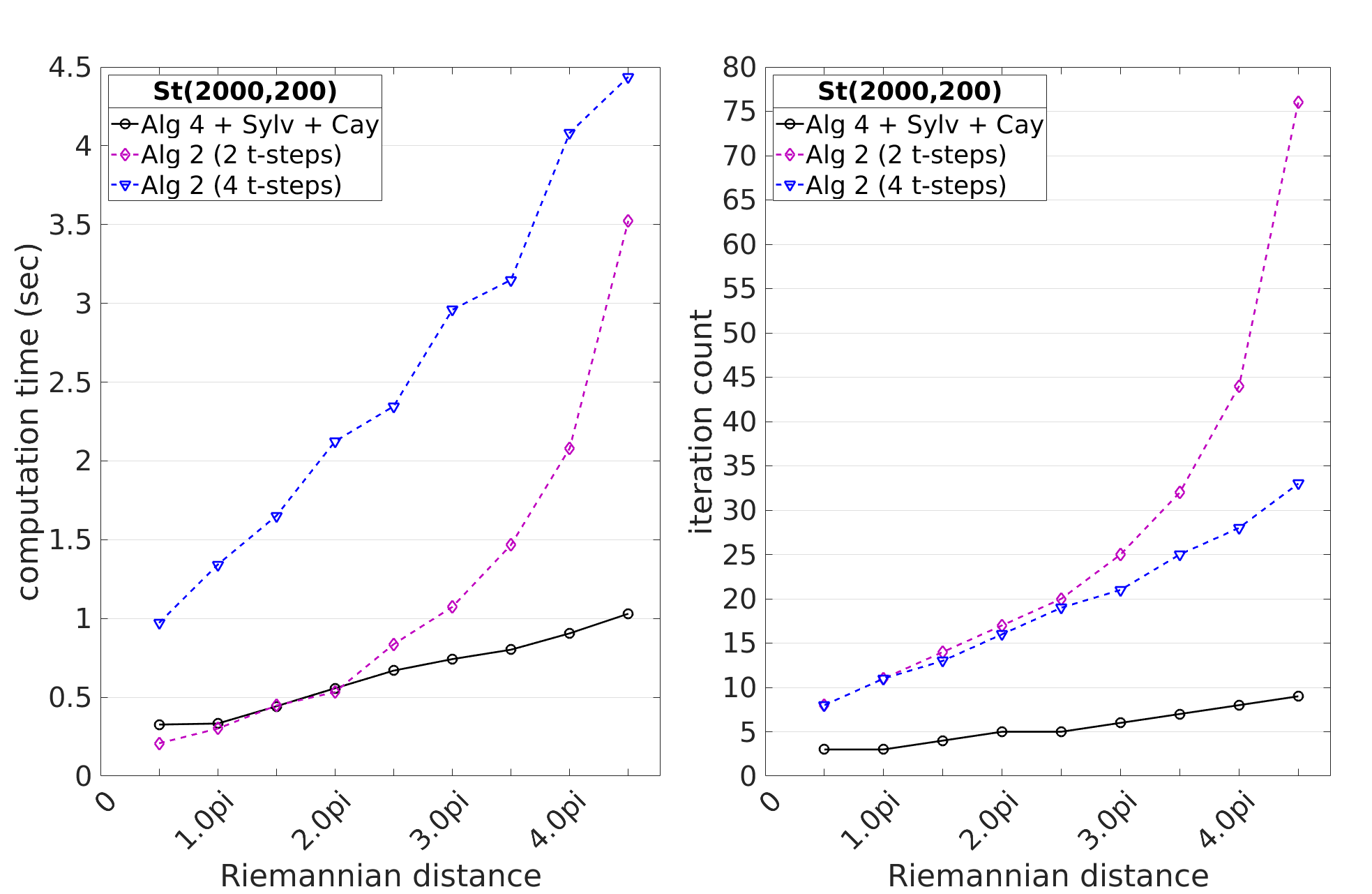}
\caption{(cf. \Cref{sec:3Experiment}) Left: Wall clock computation time in seconds versus $d=\dist(U,\widetilde U)$. 
Right: Iteration count versus $d=\dist(U,\widetilde U)$.
The trial parameter set is $\{d=(0.5+0.5j)\cdot \pi|\hspace{0.1cm} j = 0,\ldots,8\}$.
The data points are $U,\widetilde U\in St(2000,200)$. Timing results are averaged over 20 runs.}
\label{fig:dist_dependence_n2000_p200_dist05--45pi}
\end{figure}

\tc{
Next, we assess the dependency of \cref{alg:p-Shooting} and \cref{alg:StLogBCH} %for computing $\Log_U^0(\widetilde U)$ 
on the distance of the input parameters $U, \widetilde U$ under the canonical metric. 
As in \Cref{sec:1Experiment}, we construct $U\in St(n,p)$ and  $\Delta_0 \in T_USt(n,p)$.
Then, we scale $\Delta(d) = \frac{d}{\|\Delta_0\|_0} \Delta_0$ and set $\widetilde U = \Exp_U^0(\Delta(d))$ so that 
by construction, $\dist(U,\widetilde U) = d$ with respect to the canonical metric ($\alpha =0$). We consider distance factors of $d=0.5\pi, 1.0\pi, \ldots, 4.5\pi$ and measure the wall clock time and the iteration count until the numerical convergence measure drops below a threshold parameter of $\tau = 10^{-11}$.
The results are displayed in \Cref{fig:dist_dependence_n2000_p200_dist05--45pi}.
It can be seen that for \cref{alg:StLogBCH}, the iteration count and wallclock time grow moderately with increasing distance.
In contrast, time and iteration count associated with \cref{alg:p-Shooting} on two time steps exhibit a strong nonlinear dependency on the distance.
For \cref{alg:p-Shooting} on four time steps, the growth in time and number of iterations is roughly linear but with a steeper slope when compared to \cref{alg:StLogBCH}.
}

%---------------------------------------------------------------------------

\begin{figure}[ht]
\centering
\includegraphics[width=1.0\textwidth]{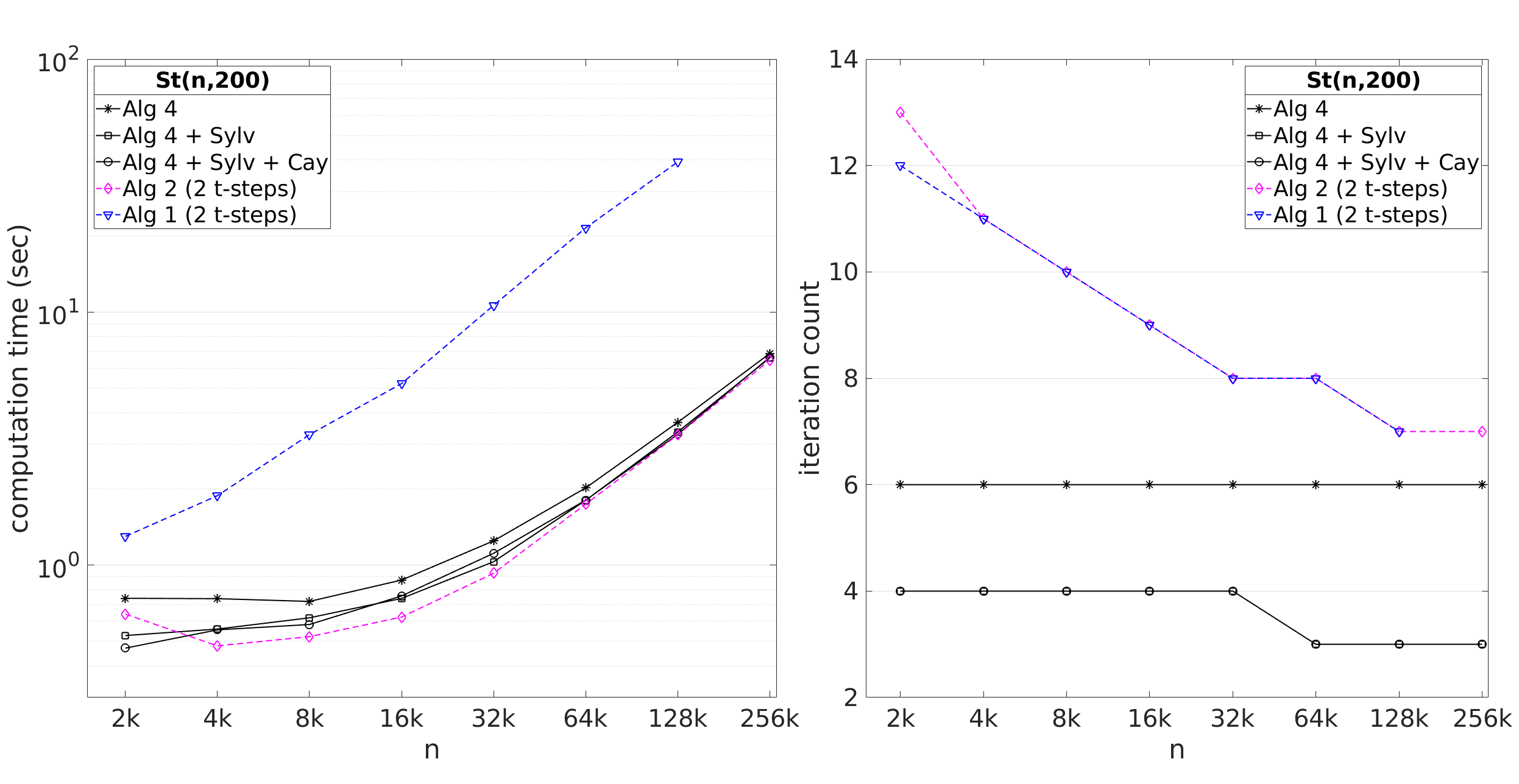}
\caption{(cf. \Cref{sec:3Experiment}) Left: Wall clock computation time in seconds versus matrix dimension $n$ on a log-log scale. 
Right: Iteration count versus matrix dimension $n$ on a log-linear scale.
The trial parameter set is $\{n=1000\cdot 2^j|\hspace{0.1cm} j = 1,\ldots,8\}$.
The data points $U,\widetilde U\in St(n,200)$ are at a canonical Riemannian distance of $\dist_\alpha(U,\widetilde U) = 1.5\pi$.
(The case of $St(256k, 200)$ was not treatable with \cref{alg:SimpleShooting} on the given laptop computer.)}
\label{fig:Stiefel_Log_n-dependencen_p200_d_1_5_pi}
\end{figure}

\tc{Lastly, we expose the dependency of \cref{alg:SimpleShooting}, \cref{alg:p-Shooting} and \cref{alg:StLogBCH} on the dimensions $n$ and $p$.
Again, we resort to the canonical metric ($\alpha =0$). We fix $p$ to a value of $p=200$ and vary $n=1000\cdot 2^j$ with $j=1,\ldots, 8$.
Pseudo-random data $U,\widetilde U \in St(n,p)$ with $\dist_0(U,\widetilde U)=1.5\pi$ is constructed as in \Cref{sec:1Experiment}. We run the various algorithms with a target accuracy of $\tau = 10^{-10}$ and measure the wall clock computation time as well as the iteration count. \cref{fig:Stiefel_Log_n-dependencen_p200_d_1_5_pi} displays the results.
The associated measurement data for $j=3,\ldots,8$ is listed in \Cref{tab:numEx_n}.
It can be observed that \cref{alg:p-Shooting} and the variants of \cref{alg:StLogBCH} perform comparably,
while \cref{alg:SimpleShooting} is roughly one order of magnitude slower.
This is expected, because the former methods address a matrix problem that scales in the dimension $p$ in the associated loop iterations, but they share matrix multiplications and a QR-decomposition of $n$-by-$p$ matrices as pre- and postprocessing steps.
}

%---------------------------------------------------------------------------
%
%

\begin{figure}[ht]
\centering
\includegraphics[width=1.0\textwidth]{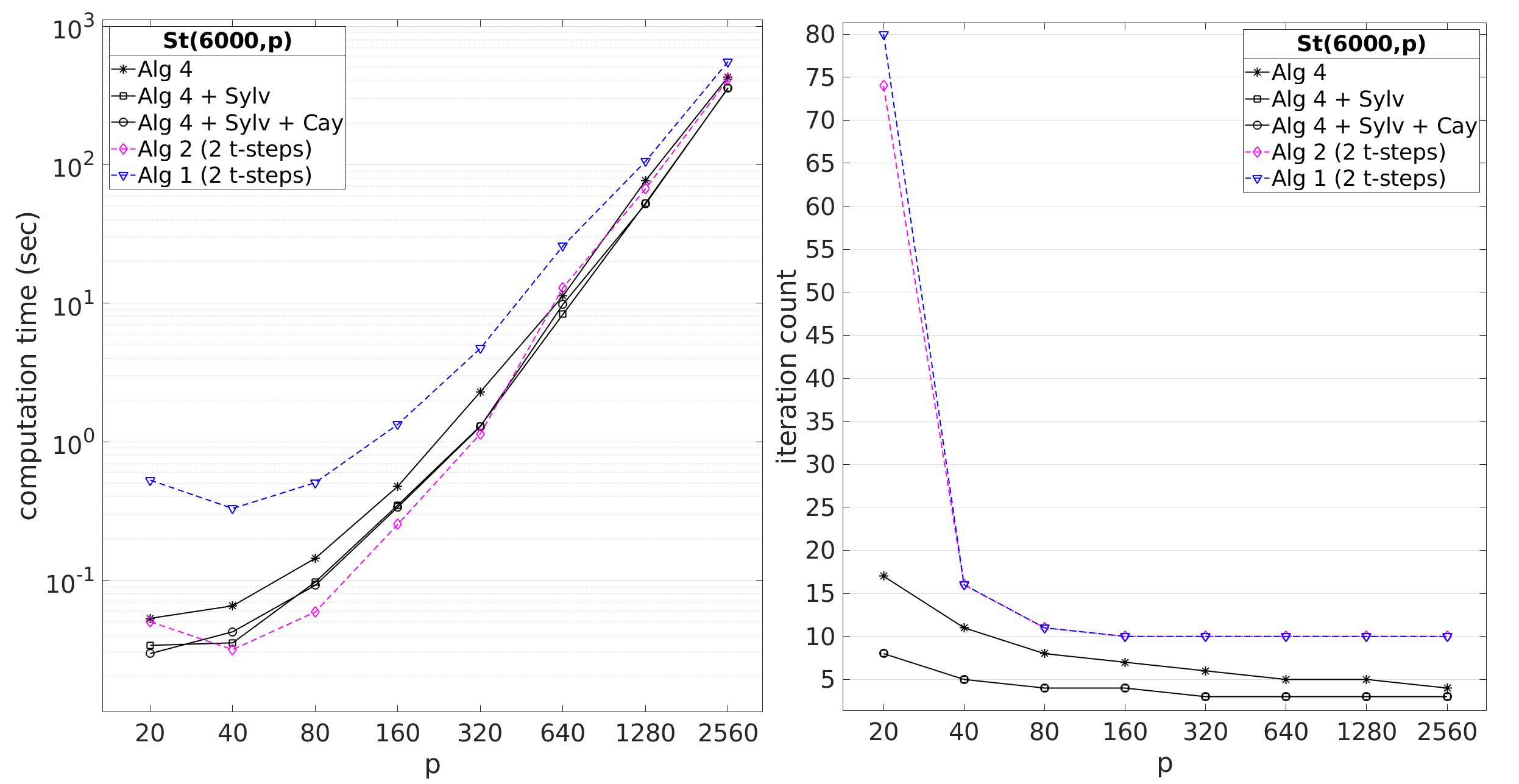}
\caption{(cf. \Cref{sec:3Experiment}) Left: Wall clock computation time in seconds versus matrix dimension $p$ on a log-log scale. 
Right: Iteration count versus matrix dimension $p$ on a log-linear scale. The trial parameter set is $\{p=10\cdot 2^j|\hspace{0.1cm} j = 1,\ldots,8\}$.
The data points $U,\widetilde U\in St(6000,p)$ are at a canonical Riemannian distance of $\dist_\alpha(U,\widetilde U) = 1.5\pi$.
}
\label{fig:Stiefel_Log_p-dependencen6000_d_1_5pi}
\end{figure}

%
%
%
%
%%%%%%******table***********table************table***********************table*******************
\begin{table}[!ht]
\begin{small}
\begin{tabular}{l|l|l|l|l|l|l}
\hline
  \rowcolor{gray!20}
Method &  $n=8k$   &$n=16k$   & $n=32k$  & $n=64k$  & $n=128k$ & $n=256k$ \\
\hline
Alg. \ref{alg:StLogBCH}                      &   0.719s &   0.873s &   1.25s &   2.02s &  3.67s  & 6.87s\\
\hline
  \rowcolor{gray!10}
Alg. \ref{alg:StLogBCH}+Sylv.                &   0.620s &   0.739s &   1.03s &   1.80s &  3.36s  & 6.59s\\
\hline
Alg.  \ref{alg:StLogBCH}+Sylv.+Cay.          &   0.582s &   0.756s &   1.11s &   1.81s &  3.29s  & 6.64s\\
\hline
  \rowcolor{gray!10}
Alg.  \ref{alg:p-Shooting} (2 t-steps)       &   0.522s &   0.623s &   0.930s&   1.74s &  3.27s  &  6.46s\\
\hline
Alg.  \ref{alg:SimpleShooting} (2 t-steps)   &   3.28s  &   5.22s  & 10.63s  &    21.5s&  39.2s  & unfeasible  
\end{tabular}
\end{small}
\caption{(cf. \Cref{sec:3Experiment}) Wallclock computation time for the results displayed in 
\cref{fig:Stiefel_Log_n-dependencen_p200_d_1_5_pi}.
}
\label{tab:numEx_n}
\end{table}
\tc{We repeat the experiment with fixed $n=6000$ and $p=10\cdot 2^j$ with $j=1,\ldots, 8$.
\cref{fig:Stiefel_Log_p-dependencen6000_d_1_5pi} displays the results.
The associated measurement data for $j=3,\ldots,8$ is listed in \Cref{tab:numEx_p}.
It can be observed that for the dimensions tested, \cref{alg:p-Shooting} is fastest until a column-dimension of $p=320$.
Beyond this point, \cref{alg:StLogBCH} with the Sylvester enhancement takes the lead.
For $p=2560$, the wallclock run time of \cref{alg:StLogBCH} is a factor of $0.82$ smaller than the runtime of \cref{alg:p-Shooting}.
}

%---------------------------------------------------------------------------
%
%

%%%%%%******table***********table************table***********************table*******************
\begin{table}[!ht]
\begin{small}
\begin{tabular}{l|l|l|l|l|l|l}
\hline
  \rowcolor{gray!20}
Method &  $p=80$   &$p=160$   & $p=320$  & $p=640$  & $p=1280$ & $p=2560$ \\
\hline
Alg. \ref{alg:StLogBCH}                      &   0.134s &   0.474s &   2.12s &  11.44s &  87.31s & 423.3s\\
\hline
  \rowcolor{gray!10}
Alg. \ref{alg:StLogBCH}+Sylv.                &   0.088s &   0.343s &   1.27s &   7.68s &  53.81s & 363.0s\\
\hline
Alg.  \ref{alg:StLogBCH}+Sylv.+Cay.          &   0.086s &   0.342s &   1.28s &   7.55s &  52.80s & 355.1s\\
\hline
  \rowcolor{gray!10}
Alg.  \ref{alg:p-Shooting} (2 t-steps)       &   0.076s &   0.226s &   1.20s &   9.87s &  61.78s & 432.3s\\
\hline
Alg.  \ref{alg:SimpleShooting} (2 t-steps)   &   0.477s &   1.336s &   4.82s &  23.28s & 112.6s & 564.0s
\end{tabular}
\end{small}
\caption{(cf. \Cref{sec:3Experiment}) Wallclock computation time for the results displayed in 
\cref{fig:Stiefel_Log_p-dependencen6000_d_1_5pi}.
}
\label{tab:numEx_p}
\end{table}
%
%
%%
%%%
%%%%
%%%%%
\subsection{The impact of the metric on interpolating matrix factorizations}
\label{sec:matInt}
% # QR or SVD
% key = 'QR'
% 
% # metric
% alpha_range = np.arange(-0.8,2.05,0.02)
%results
% Tang RBF: min L2 entry: val  0.019290355763450954 index 54
% alpha val= 0.28
% Cubic Spline: min L2 entry: val  0.0161644447357338 index 53
% alpha val= 0.26
In this subsection, we consider the practical problems of interpolating the QR-decom\-po\-sition and 
the SVD of time-dependent matrix curves. Let $n>m>p$. For a compact QR-decomposition $Y=QR$ of a rectangular matrix $Y\in \R^{n\times p}$, it holds $Q\in St(n,p)$. If the SVD is used to produce the best rank-$p$ approximation of a given matrix $\R^{n\times m}\ni Y\approx U_p\Sigma_p V^T_p$, then $U_p\in St(n,p), V_p\in St(m,p)$. 
Hence, in both applications, interpolation of matrix sample data on the Stiefel manifold has to be considered. 
We will investigate how the choice of the metric affects the resulting interpolant.

The standard approach to interpolating manifold-valued data is
(1) to select a base point, 
(2) to apply the Riemannian logarithm to map the data set to the tangent space at the chosen base point, 
(3) to perform interpolation in the tangent space, 
(4) to map the results back to the manifold via the Riemannian exponential.

First, we reproduce the example from \cite[Section 5.2]{ZimmermannHermite_2020}
and consider a cubic matrix polynomial
\[ 
 t\mapsto Y(t) = Y_0 + t Y_1 + t^2Y_2 + t^3Y_3, \quad Y_k \in \R^{n\times p}, \quad  n=500, p=10.
\]
The matrices $Y_k$ are produced as random matrices with entries uniformly sampled from $[0,1]$ for $Y_0$,
entries uniformly sampled from $[0,0.5]$ for $Y_1,Y_2$ and from $[0,0.2]$ for $Y_3$.
Then, $Y(t)$ is sampled at $5$ equidistant samples $t_i \in \{-1.1,  -0.55,  0.0,    0.55,  1.1 \}$.
At each sample point $t_i$, the Q-factor $Q(t_i)\in St(n,p)$ of the QR-decomposition is computed. 
We employ radial basis function (RBF) interpolation in the tangent space with the cubic RBF 
and three-point cubic spline interpolation, for details, see \cite{ZimmermannHermite_2020}.
For mapping the sample data back and forth between the Stiefel manifold and its tangent space as sketched in \Cref{fig:ExpLogStiefel}, we use the Riemannian exponential and logarithm under the $\alpha$-metric.
We start at $\alpha = -0.8$ and proceed until $\alpha = 2.04$ in steps of $\delta\alpha =0.02$.

\begin{figure}[ht]
\centering
\includegraphics[width=0.8\textwidth]{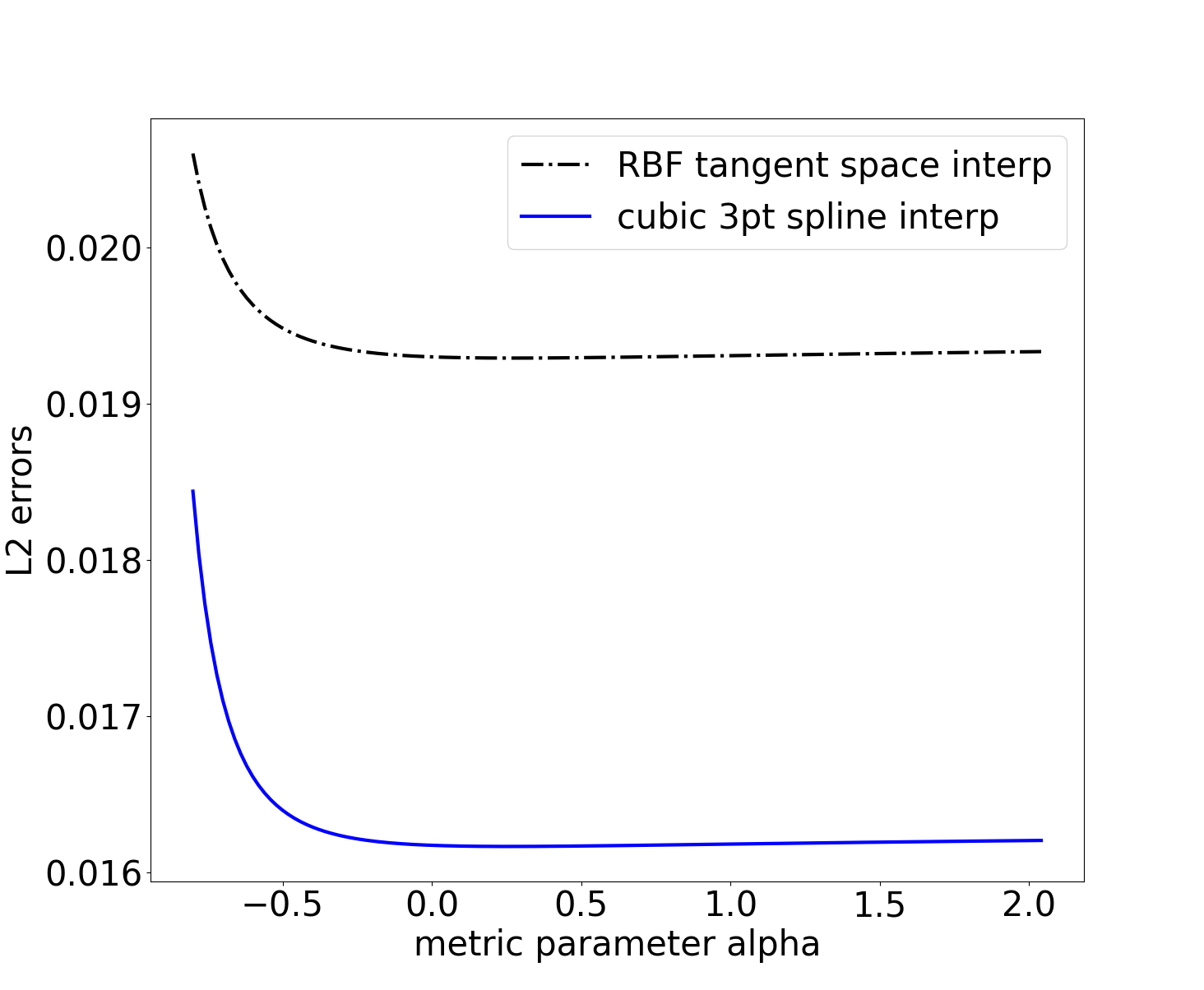}
\caption{(cf. \Cref{sec:matInt})Total $L_2$ errors versus $\alpha$ associated with interpolating the $Q$-factor of a time-dependent curve of QR-factorizations under varying the $\alpha$-metric. The tested $\alpha$-range is $[-0.8, 2.04]$ on steps of size $\delta\alpha = 0.02$.
}
\label{fig:matInt_QR}
\end{figure}

For each value of $\alpha$, the relative interpolation errors are computed at $101$ equidistant instants $t_j=-1.1+j\delta t\in [-1.1,1.1]$,
$\delta t = 0.022$, $j=0,1,\ldots,100$
in the matrix Frobenius norm
\[
    e_\alpha(t_j) := \frac{\|Q^*_\alpha(t_j) - Q(t_j)\|_F}{\|Q(t_j)\|_F}.
\]
Here, $Q^*_\alpha(t_j)$ denotes the manifold interpolant under the $\alpha$-metric and $Q(t_j)$ is the reference solution.
The associated discrete $L_2$-norm of the error is computed as $\text{err}_{L_2}(\alpha) = \sqrt{\delta t\sum_{j=0}^{100} e_\alpha(t_j)^2}$.
The error graphs $\alpha\mapsto \text{err}_{L_2}(\alpha)$ for both interpolation methods under consideration are displayed in \Cref{fig:matInt_QR}.
%---------------------------------------------------------------------------
The discrete arrays that underlie the graphs feature both a global minimum at a similar location.
For RBF interpolation, the total $L_2$ error is lowest at $\alpha = 0.28$.
For three-point cubic spline interpolation, the global minimum is at $\alpha = 0.26$.
It can also be seen that beyond $\alpha \geq -0.1$, the impact of the metric on the error is rather negligible.
It should also be emphasized, that the largest and the smallest tested errors differ only by a small absolute amount.

As a second test case, we rework the example from \cite[Section 5.3]{ZimmermannHermite_2020} and construct a nonlinear matrix function with fixed low rank. We start with a cubic matrix polynomial
\[ 
 Y(t) = Y_0 + t Y_1 + t^2Y_2 + t^3Y_3, \quad Y_i \in \R^{n\times r}, n=10\hspace{0.07cm}000, r=10
\]
with random matrices $Y_k$ with entries uniformly sampled from $[0,1]$ for $Y_0$
and from $[0,0.5]$ for $Y_1,Y_2, Y_3$.
Then, a second matrix polynomial is considered
\[ 
 Z(t) = Z_0 + t Z_1 + t^2Z_2, \quad Z_i \in \R^{r\times m}, r=10, m=300.
\]
Here, the entries of $Z_0$ are sampled uniformly from $[0,1]$
while the entries of $Z_1,Z_2$ are sampled uniformly from $[0,0.5]$.
The nonlinear low-rank matrix function is set to be
\[
 W(t) = Y(t) Z(t) \in \R^{n\times m}.
\]
We will conduct cubic Hermite interpolation. To this end, the low-rank SVD
\[
 W(t) = U_r(t)\Sigma_r(t)V_r(t)^T, \quad U_r(t)\in St(n,r), \hspace{0.2cm} V_r(t)\in St(m,r),\hspace{0.2cm}  \Sigma \in \R^{r\times r}
\]
and the associated matrix derivatives are sampled at three Chebychev nodes $t_0\approx0.0603, t_1 = 0.45, t_2 = 0.8397$ in the interval $[0.0, 0.9]$. 
The Hermite sample data set is
\[
 U_r(t_i), \dot U_r(t_i), \quad V_r(t_i), \dot V_r(t_i),\quad \Sigma_r(t_i), \dot \Sigma_r(t_i),\quad i=0,1,2,
\]
see \cite[Section 5.3]{ZimmermannHermite_2020} for details.

We conduct Hermite interpolation under the $\alpha$-metric, starting from $\alpha =-0.75$ and proceeding up to $\alpha = 1.5$ in steps of $\delta\alpha =0.05$. (For $\alpha < -0.75$ and $\alpha > 1.5$, convergence issues occurred when mapping the sample data set to the tangent space with the Riemannian logarithm.)

For each value of $\alpha$, the relative interpolation errors are computed at $100$ equidistant instants $t_j\in [t_0,t_2]$
in the matrix Frobenius norm as
$$
e_\alpha(t_j) = \frac{\|U^*(t_j)\Sigma^*(t_j)(V^*(t_j))^T - W(t_j)\|_F}{\|W(t_j)\|_F},
$$ 
where $U^*(t_j)$, $\Sigma^*(t_j)$, $V^*(t_j)$ are the interpolants of the 
matrix factors of the low-rank SVD of $W(t_j)$ and $W(t_j)$ is the reference solution.
%---------------------------------------------------------------------------
\begin{figure}[ht]
\centering
\includegraphics[width=0.8\textwidth]{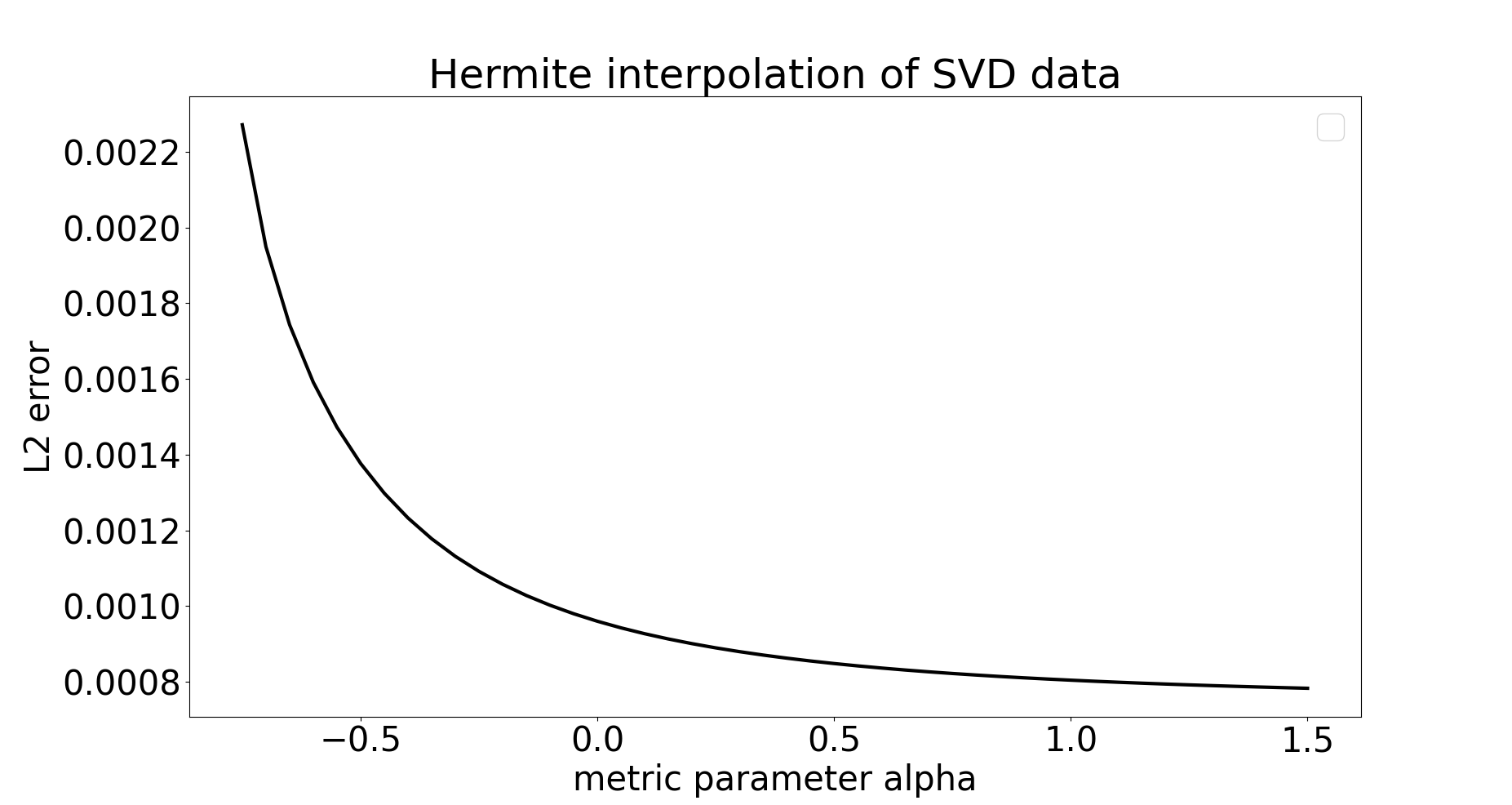}
\caption{(cf. \Cref{sec:matInt}) Total $L_2$ errors versus $\alpha$ associated with Hermite-interpolating a time-dependent curve of low-rank SVDs under varying the $\alpha$-metric. The tested $\alpha$-range is $[-0.75, 1.5]$ on steps of size $\delta\alpha = 0.05$.
}
\label{fig:matInt_SVD}
\end{figure}
The associated discrete $L_2$-norm of the error is computed as $\text{err}_{L_2}(\alpha) = \sqrt{\delta t\sum_{j=0}^{99} e_\alpha(t_j)^2}$.
The error graph $\alpha\mapsto \text{err}_{L_2}(\alpha)$ is displayed in \Cref{fig:matInt_SVD}.
It can be seen that in the tested $\alpha$-range, the total $L_2$ error is largest at the left bound $\alpha=-0.75$ 
and decreases monotonically with increasing $\alpha$. The slope flattens considerably for larger values of $\alpha$.
%
%%
%%%
%%%%
%%%%%
%%%%%%
%------------------------------------------------------
\section{Discussions}
\label{sec:conclusions}

For the canonical metric, the numerical experiments identify the Sylvester-en\-hanced version of \Cref{alg:StLogBCH}  as the method of choice:
The runtime is comparable to \Cref{alg:p-Shooting} on two time steps for the smaller test cases and considerably shorter for the test cases in higher dimensions. Moreover, it is more robust with respect to increasing the distance of the input points, see \Cref{fig:dist_dependence_n2000_p200_dist05--45pi}.

For all other $\alpha$-metrics, the $p$-shooting method \Cref{alg:p-Shooting} on two time steps performs best in terms of the runtime. Yet, data on Stiefel manifolds of smaller dimension or data points that are further apart may necessitate in doing more time steps in the inner loops of \Cref{alg:p-Shooting}.
Even though for all test cases considered in this work, it was sufficient to go up to $4$ times steps, \Cref{alg:p-Shooting} is a local method and it cannot be expected that increasing the number of inner time steps will make the method  converge in all cases.
The Euclidean metric ($\alpha = -\frac12$) has a special position in the sense that the iteration count of the $p$-shooting method is lowest for the Stiefel manifold under this choice of metric. 

The choice of metric is problem-dependent.
We included an experiment on interpolating curves of orthogonal matrix factorizations.
When casting this into an interpolation problem on the Stiefel manifold, a metric has to be selected.
While the exact curve is completely independent from the selected metric, the interpolation error is not.
Hence, if the computational resources allow for it, a parametric study may be beneficial when selecting the metric for a certain application. Otherwise, the canonical metric presents itself as a natural general purpose tool.
%
%
%%
%%%
%%%%
%%%%%
%%%%%%
%%%%%%%
%%%%%%%%
%%%%%%%%%
%%%%%%%%%%
% 
% Curvature Estimate for Stiefel: EUCLIDEAN METRIC: \cite[Prop. 3.4]{Chakraborty2017}
% Proposition 3.4. Let 
% $X\in St(n,r), \Delta,\widetilde{\Delta}\in \mathcal{H}_X\cong T_{[X]}Gr(n,r)$, $\sigma = span(\Delta,\widetilde{\Delta})$, then, $0 \leq \kappa_X(\sigma) \leq 2$.
% CANONICAL METRIC:
% RENTMEESTERS, p. 99, 100
%
\section*{Acknowledgements}
The authors would like to thank Marco Sutti for sharing the Matlab code assoicated with \cite{Sutti_PhD_2020}.
The second author's work has been supported in part by the German Federal Ministry of Education and Research (BMBF-Projekt 05M20WWA: Verbundprojekt 05M2020 - DyCA).

Independent from this work and concurrent with our original arXiv submission, the related preprint \cite{nguyen2021closedform} appeared, which takes a similar starting point to tackle the geodesic endpoint problem but then pursues an optimization approach.
%
%
%
%
%

%
%
%%
%%%
%%%%
%%%%%
%%%%%%
\appendix
\section{Convergence plots associated with \texorpdfstring{\Cref{sec:1Experiment}}{Subsection 4.1}}
\label{supp:1Experiment}
\begin{figure}[ht]
\centering
\includegraphics[width=0.9\textwidth]{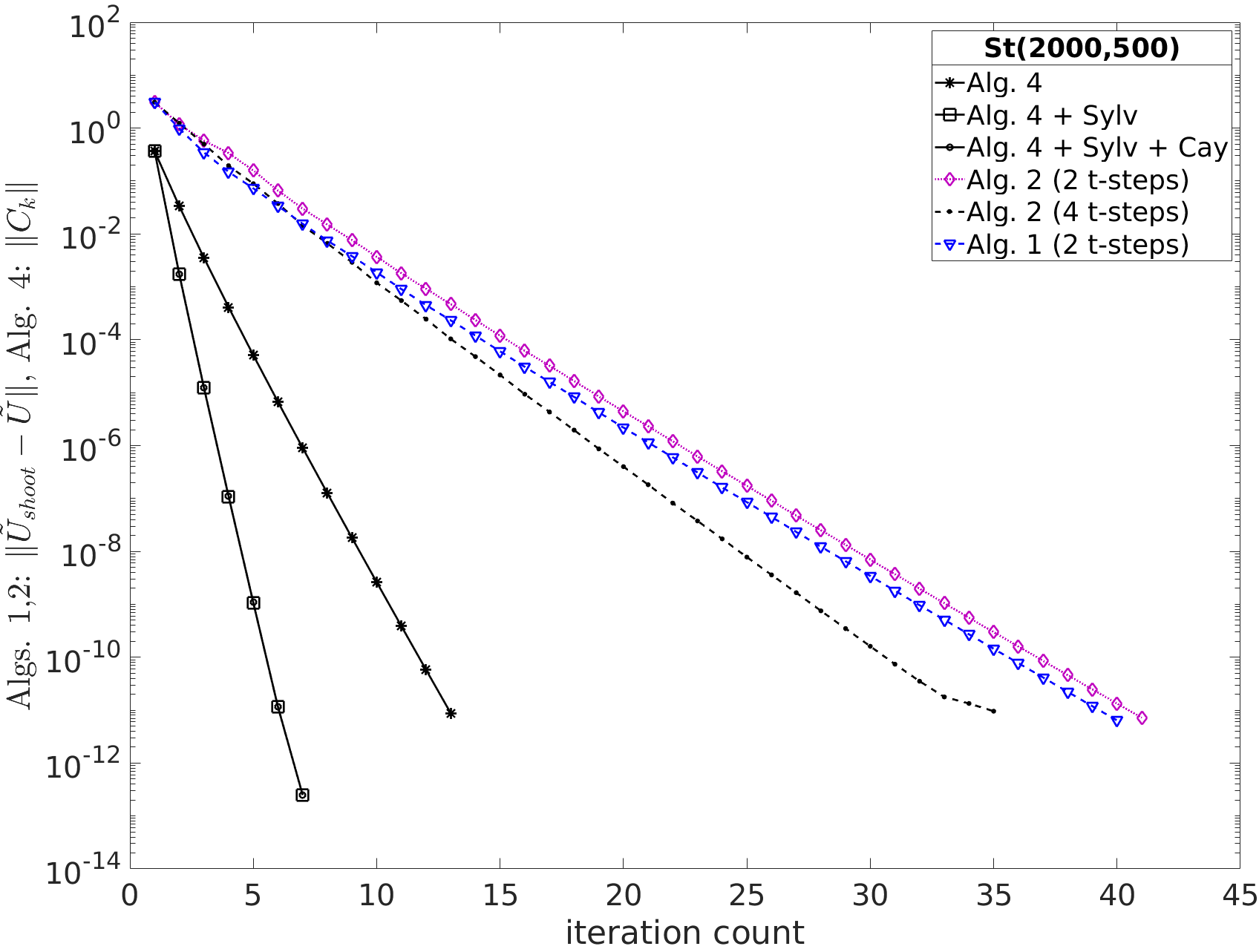}
\caption{(cf. \Cref{sec:1Experiment}) Convergence history of the various log-algorithms for computing $\Log_U(\widetilde U)$ for the canonical metric. The graphs show one run of the test case $St(2000, 500)$.
The input data is at a Riemannian distance of $\dist(U,\widetilde U) = 5\pi$.
}
\label{supp:fig:St_canoni_2000_500_5pi}
\end{figure}

\begin{figure}[ht]
\centering
\includegraphics[width=0.9\textwidth]{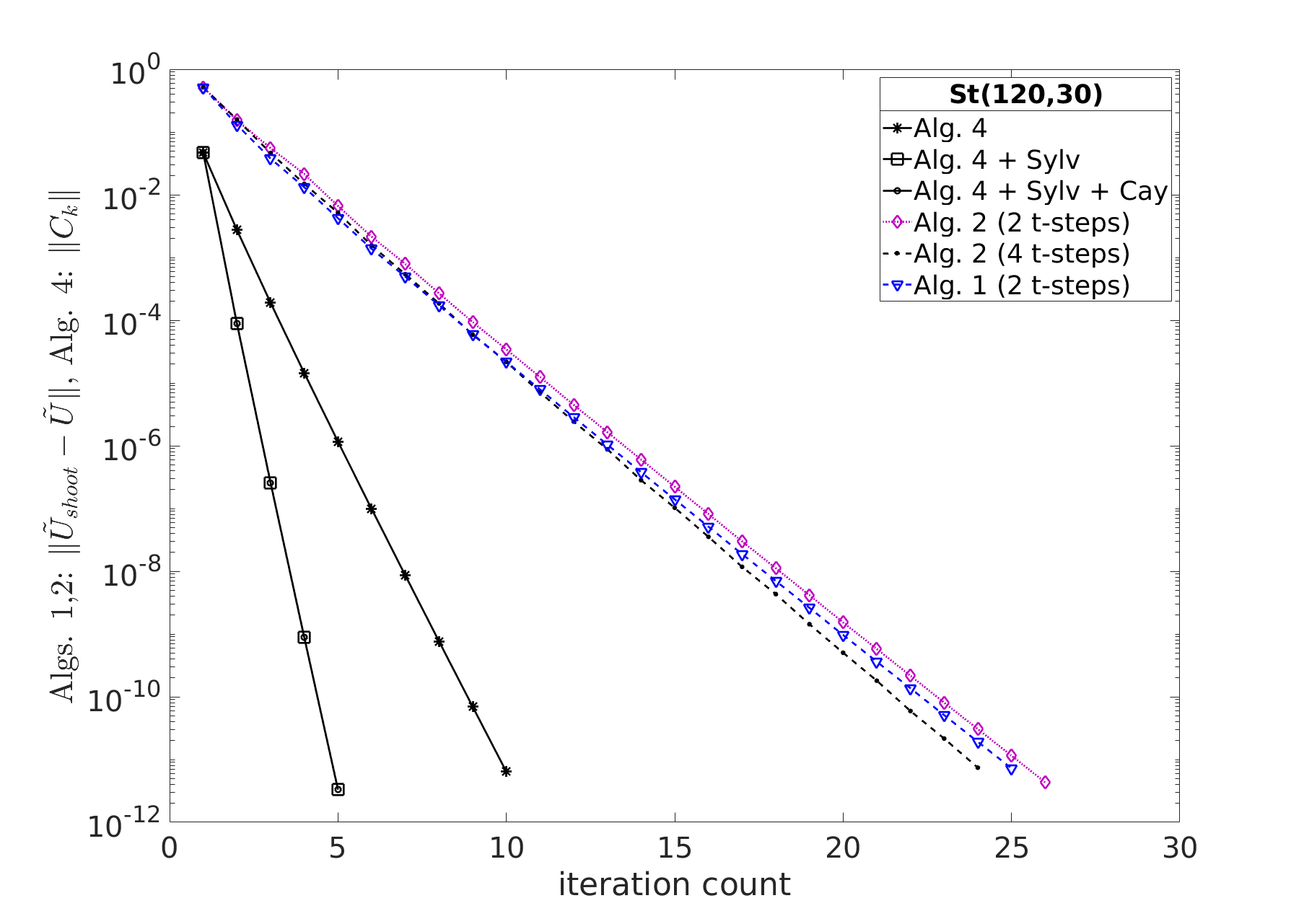}
\caption{(cf. \Cref{sec:1Experiment}) Convergence history of the various log-algorithms for computing $\Log_U(\widetilde U)$ under the canonical metric. The graphs show one run of the test case $St(120, 30)$.
The input data is at a Riemannian distance of $\dist(U,\widetilde U) = \pi$.
}
\label{supp:fig:St_canoni_120_30_1pi}
\end{figure}

\begin{figure}[ht]
\centering
\includegraphics[width=0.9\textwidth]{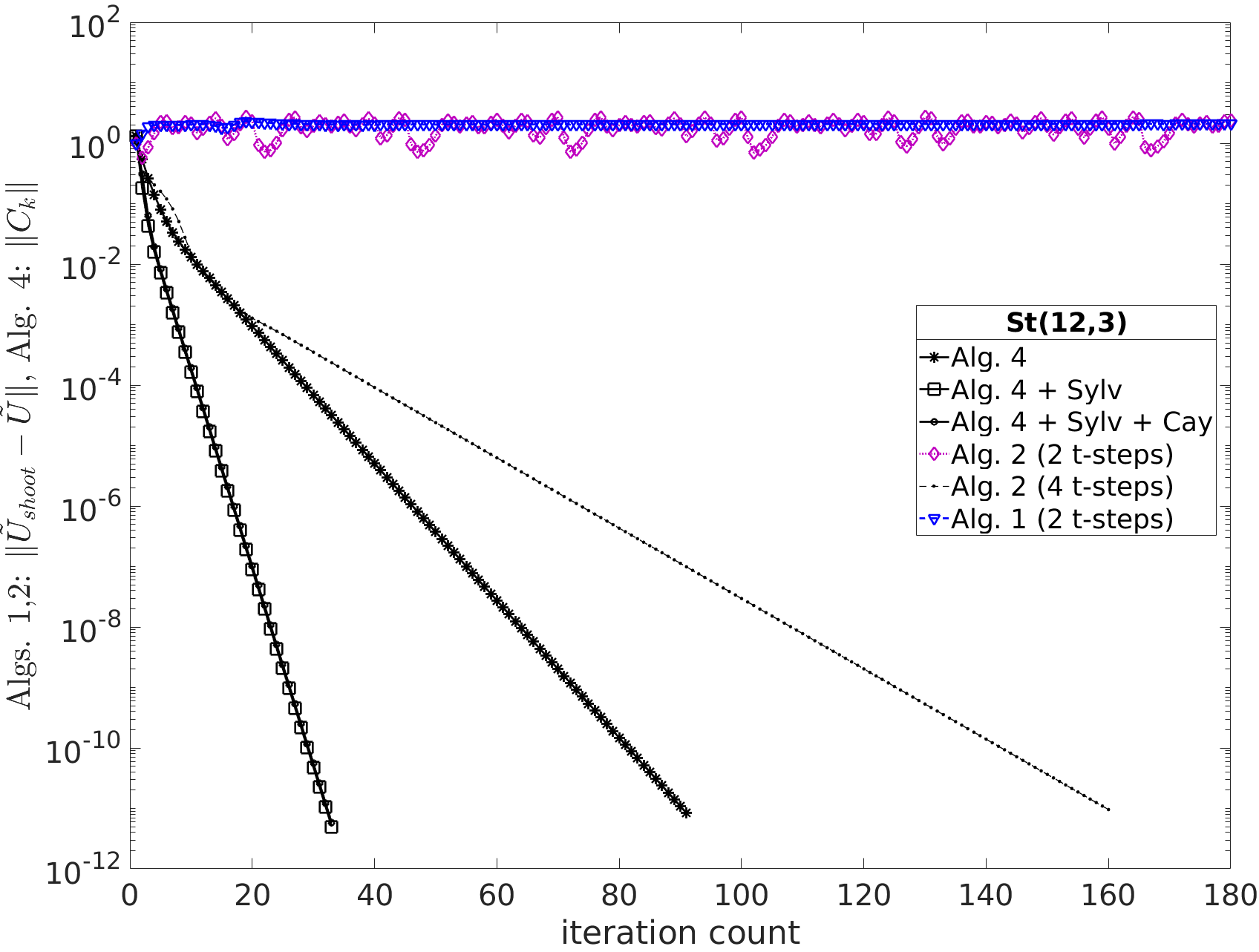}
\caption{(cf. \Cref{sec:1Experiment}) Convergence history of the various log-algorithms for computing $\Log_U(\widetilde U)$ under the canonical metric. The graphs show one run of the test case $St(12, 3)$.
The input data is at a Riemannian distance of $\dist(U,\widetilde U) = 0.95\pi$.
}
\label{supp:fig:St_canoni_12_3_1pi}
\end{figure}

Figures \ref{supp:fig:St_canoni_2000_500_5pi}, \ref{supp:fig:St_canoni_120_30_1pi} and \ref{supp:fig:St_canoni_12_3_1pi}
show the convergence history of a single run associated with  Test Case 1, Test Case 2 and Test Case 3 of \Cref{tab:numEx1}, respectively. Note that for \Cref{alg:StLogBCH}, the convergence measure at iteration $k$ is the norm of the lower diagonal block $C_k$  in \eqref{eq:basicStLog},
while it is the error between the target $\widetilde U$ and the shooting point $\widetilde U_{\text{shoot}}$ at iteration $k$ for \Cref{alg:p-Shooting} and \Cref{alg:SimpleShooting}.

%---------------------------------------------------------------------------
%---------------------------------------------------------------------------

%St_canoni_120_30_1pi.png
%%
%%
% 
% %---------------------------------------------------------------------------
% \begin{figure}[ht]
% \centering
% \includegraphics[width=1.0\textwidth]{LogPerformance_n100000_p400_d5pi.png}
% \caption{(cf. \Cref{sec:1Experiment}) Convergence history of the various log-algorithms for computing $\Log_U(\widetilde U)$ for the canonical metric. The graphs show one run of the test case $St(100 000, 400)$.
% The input data is at a Riemannian distance of $\dist(U,\widetilde U) = 5\pi$.
% }
% \label{fig:LogPerformance_n100000_p400_d5pi}
% \end{figure}
% %---------------------------------------------------------------------------
%
%

\section{Convergence plots associated with \texorpdfstring{\Cref{sec:2Experiment}}{Subsection 4.2}}
\label{supp:2Experiment}
\begin{figure}[hb]
\centering
\includegraphics[width=0.9\textwidth]{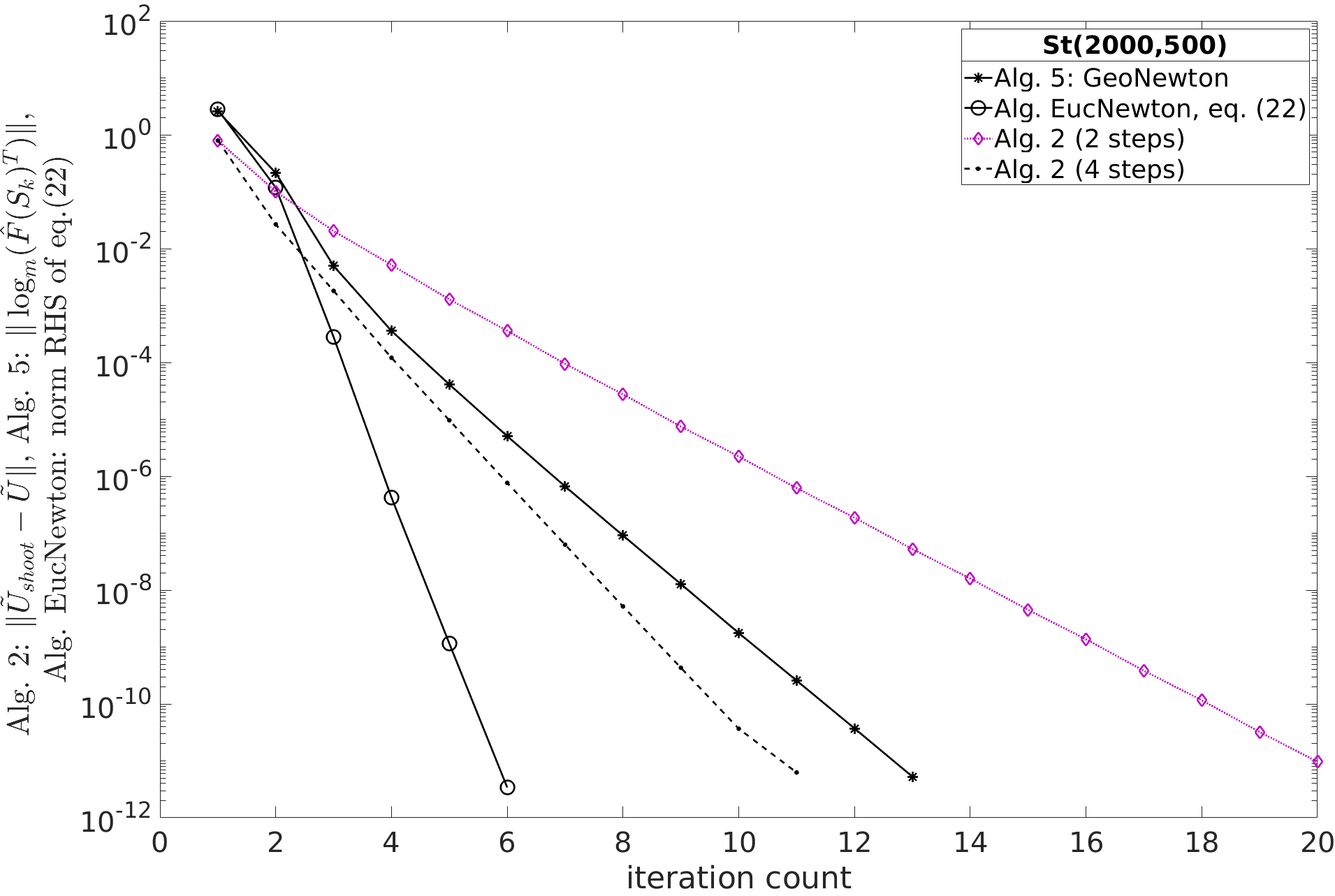}
\caption{(Cf. \Cref{sec:2Experiment}) Convergence history of the various log-algorithms for computing $\Log_U(\widetilde U)$ under the Euclidean metric. The graphs show one run of the test case $St(2 000, 500)$.
The input data is at a Riemannian distance of $\dist(U,\widetilde U) = 5\pi$.
}
\label{fig:St_Euclid_2000_500_5pi}
\end{figure}

\begin{figure}[hb]
\centering
\includegraphics[width=0.9\textwidth]{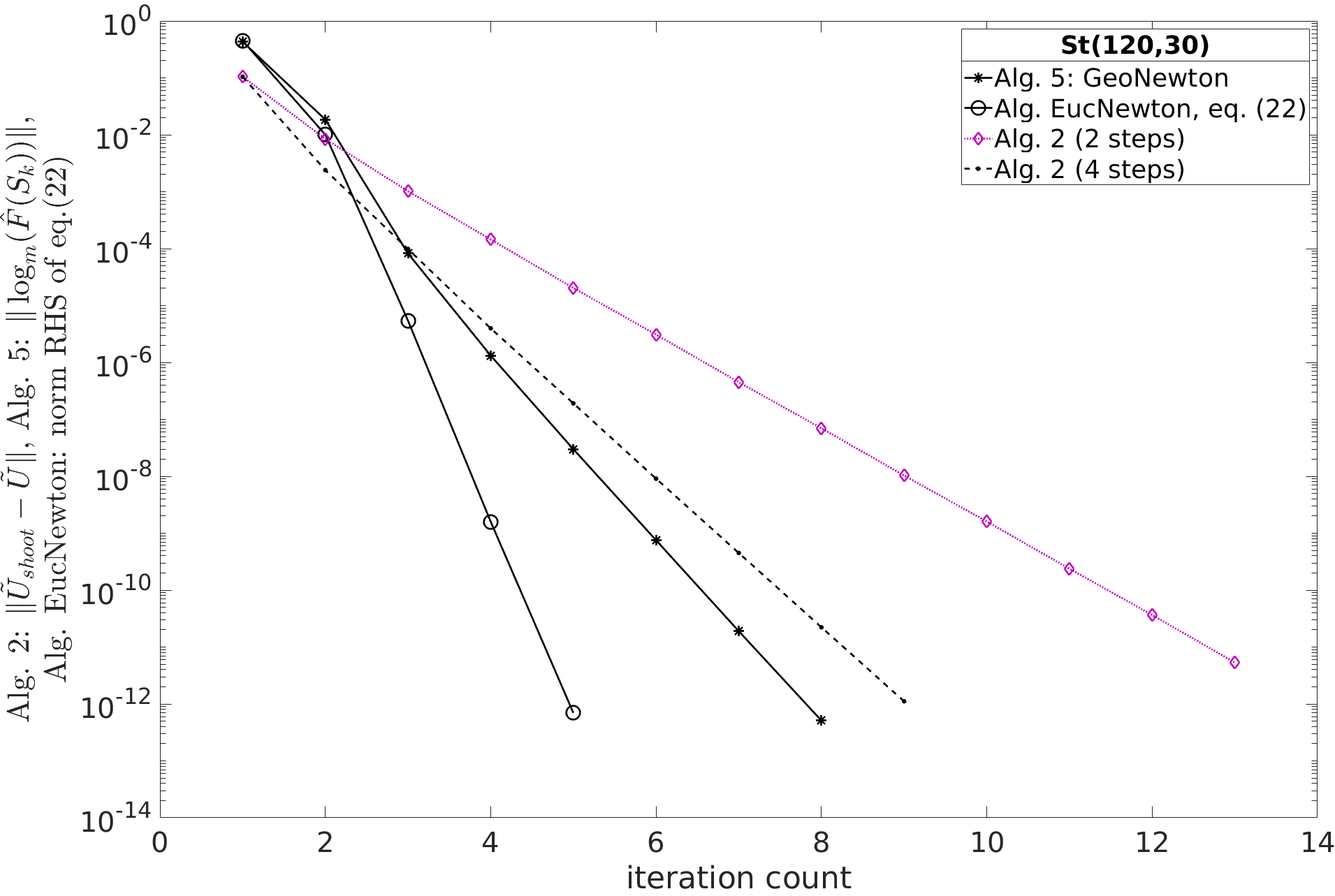}
\caption{(Cf. \Cref{sec:2Experiment}) Convergence history of the various log-algorithms for computing $\Log_U(\widetilde U)$ under the Euclidean metric. The graphs show one run of the test case $St(120, 30)$.
The input data is at a Riemannian distance of $\dist(U,\widetilde U) = \pi$.
}
\label{fig:St_Euclid_120_30_1pi}
\end{figure}

Figures \ref{fig:St_Euclid_2000_500_5pi} and \ref{fig:St_Euclid_120_30_1pi}
show the convergence history of a single run associated with  Test Case 1 and Test Case 2 of \Cref{tab:numEx2}, respectively. Note that the convergence measure for \Cref{alg:p-Shooting} and \Cref{alg:SimpleShooting}
is the error between the target $\widetilde U$ and the shooting point $\widetilde U_{\text{shoot}}$ at iteration $k$,
while it is $\| \log_m(\widehat F(S_k)^T)\|$ for \Cref{alg:GeoNewt} and the norm of the right hand side of \eqref{eq:EucNewton4Log}
for Algorithm `EucNewton'.

%---------------------------------------------------------------------------
%---------------------------------------------------------------------------
\section{Additional experiments associated with \texorpdfstring{\Cref{sec:3Experiment}}{Subsection 4.3}}
\label{supp:3Experiment}

\begin{figure}[ht]
\centering
\includegraphics[width=1.0\textwidth]{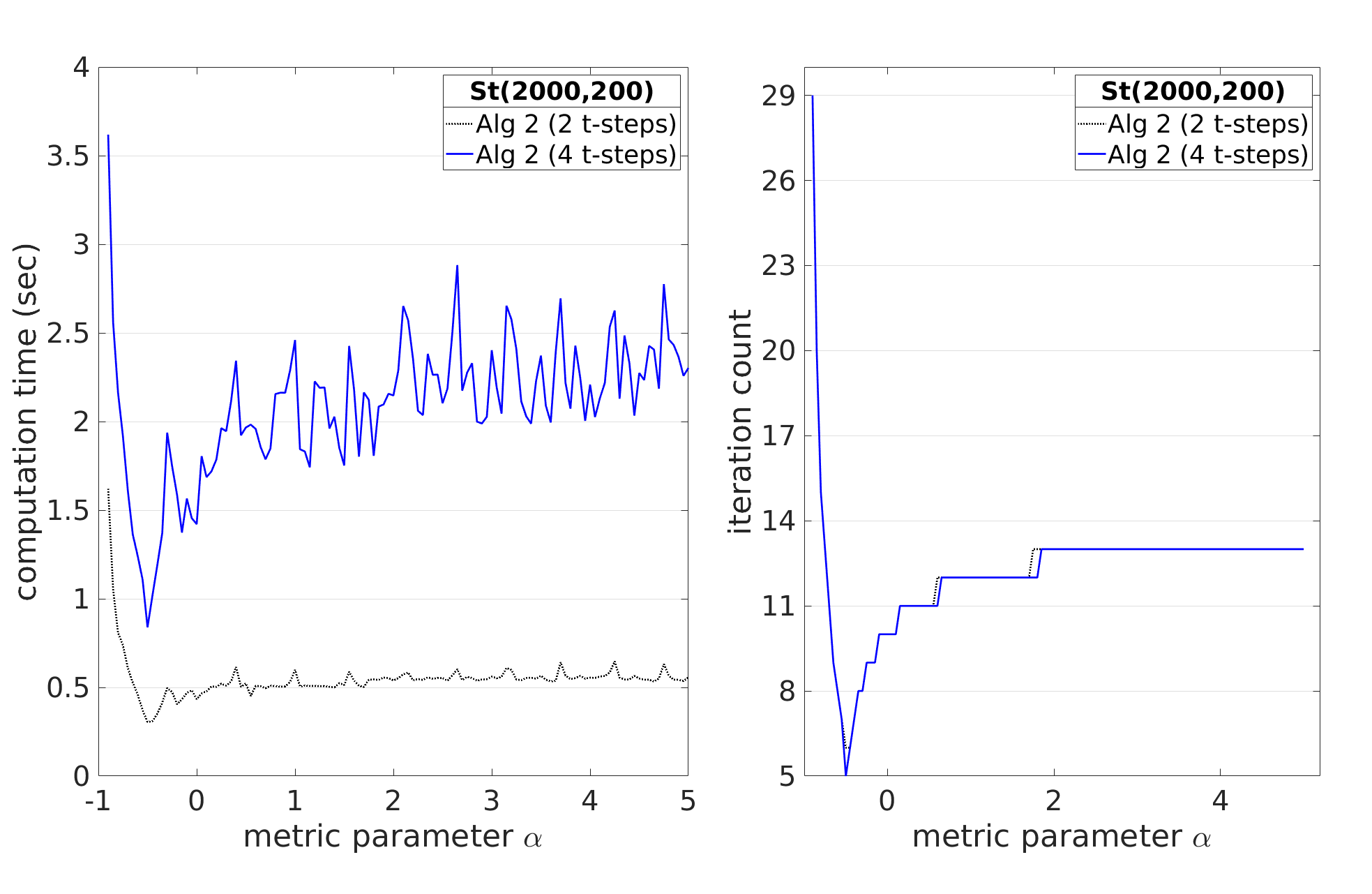}
\caption{(cf. \Cref{sec:3Experiment}) Wall clock computation time in seconds versus $\alpha$ (left) and iteration count until convergence versus $\alpha$ (right) for computing $\Log_U^\alpha(\widetilde U)$ with \Cref{alg:p-Shooting}. The trial parameter set is $\{\alpha=-0.9 + 0.05j|\hspace{0.1cm} j = 0,\ldots,118\}$.
The data points $U,\widetilde U\in St(2000,200)$ are at a Riemannian $\alpha$-distance of $\dist_\alpha(U,\widetilde U) = 0.8\pi$. Timing results are averaged over 10 runs.
}
\label{fig:alpha_dependence_n2000_p200_d08pi_alpha_-09--5}
\end{figure}

We repeat the first experiment of \Cref{sec:3Experiment} for pseudo-random data
$U,\widetilde U\in St(2000,200)$.
As a distance factor, we use $d=0.8\pi$. 
For each value $\alpha\in \{-0.9 + 0.05j| j = 0,\ldots,118\}$, we apply \cref{alg:p-Shooting} to compute 
$\Log_U^\alpha(\widetilde U)$ under the $\alpha$-metric and record the wall clock computation time as well as the iteration count.
The results are displayed in \Cref{fig:alpha_dependence_n2000_p200_d08pi_alpha_-09--5}.
%
%
%---------------------------------------------------------------------------
\bibliographystyle{plain}
\bibliography{SurrogateModelling}
\end{document}